\documentclass[pdflatex,sn-mathphys-num]{sn-jnl}

\usepackage{graphicx}%
\usepackage{multirow}%
\usepackage{amsmath,amssymb,amsfonts}%
\usepackage{amsthm}%
\usepackage{mathrsfs}%
\usepackage[title]{appendix}%
\usepackage{xcolor}%
\usepackage{textcomp}%
\usepackage{manyfoot}%
\usepackage{booktabs}%
\usepackage{algorithm}%
\usepackage{algorithmicx}%
\usepackage{algpseudocode}%
\usepackage{listings}%

\renewcommand{\normalsize}{\fontsize{12pt}{15pt}\selectfont}  

\def\dbar{{\mkern3mu\mathchar'26\mkern-12mu d}}

\usepackage{wrapfig}    

\usepackage{natbib}

\usepackage{hyperref}

\usepackage{amsfonts,amsmath,amsthm}

\usepackage{mathtools}
\newcommand\Mycomb[2][^n]{\prescript{#1\mkern-0.5mu}{}C_{#2}}

\usepackage[thinc]{esdiff}

\usepackage{leftidx}

\usepackage{longtable}   

\usepackage{physics}

\usepackage{amssymb}

\usepackage{comment}

\usepackage{lscape}

\usepackage{xcolor}

\usepackage{rotating}  

\usepackage{booktabs}
\usepackage{lscape}

\usepackage{tikz}
\usepackage{caption}

\title[Entropy measures and their applications: A comprehensive review]{Entropy measures and their applications: A comprehensive review}

\author[1]{\fnm{Naveen} \sur{Kumar}}\email{kumar.248@iitj.ac.in}
\author[2]{\fnm{Ambesh} \sur{Dixit}}\email{ambesh@iitj.ac.in}
\author[3]{\fnm{Vivek} \sur{Vijay}}\email{vivek@iitj.ac.in}

\affil[1,3]{\orgdiv{Department of Mathematics}, \orgname{IIT Jodhpur}, \orgaddress{\street{Jodhpur}, \postcode{342030}, \state{Rajasthan}, \country{India}}}
\affil[2]{\orgdiv{Department of Physics}, \orgname{IIT Jodhpur}, \orgaddress{\street{Jodhpur}, \postcode{342030}, \state{Rajasthan}, \country{India}}}

\theoremstyle{thmstyleone}%
\theoremstyle{thmstyletwo}%

\theoremstyle{thmstylethree}%

\begin{document}









\abstract{Entropy has emerged as a dynamic, interdisciplinary, and widely accepted quantitative measure of uncertainty across different disciplines. A unified understanding of entropy measures, supported by a detailed review of their theoretical foundations and practical applications, is crucial to advance research across disciplines. 
This review article provides motivation, fundamental properties, and constraints of various entropy measures. These measures are categorized with time evolution ranging from Shannon entropy generalizations, distribution function theory, fuzzy theory, fractional calculus to graph theory, all explained in a simplified and accessible manner. These entropy measures are selected on the basis of their usability, with descriptions arranged chronologically. We have further discussed the applicability of these measures across different domains, including thermodynamics, communication theory, financial engineering, categorical data, artificial intelligence, signal processing, and chemical and biological systems, highlighting their multifaceted roles. A number of examples are included to demonstrate the prominence of specific measures in terms of their applicability. The article also focuses on entropy-based applications in different disciplines, emphasizing openly accessible resources. Furthermore, this article emphasizes the applicability of various entropy measures in the field of finance. The article may provide a good insight to the researchers and experts working to quantify uncertainties, along with potential future directions.}



\maketitle



\section{Introduction}
In every field of science, researchers endeavour to understand the complex behaviours of systems, often confronted with the challenge of measuring uncertainty and complexity. Entropy has emerged as a remarkable concept to tackle this challenge. Originally grounded in thermodynamics, entropy has grown into a flexible and widely used measure across many disciplines. Its meaning and importance have evolved over time, as Camacho et al.\citep{flores2015conceptintro1} eloquently explore, revealing many layers and interpretations that entropy has acquired throughout its rich history.

\par The journey of entropy began in $1824$ with Carnot’s theory of the heat engine revealing that the efficiency of such engines depends on heat transformations. This idea was further developed by Clausius\citep{cropper1986rudolfintro2} in $1854$, who formally introduced the concept of entropy to quantify the irreversible unused heat transfer in Carnot engines. In $1872$, Boltzmann\citep{chakrabarti2000boltzmann} bridged the microscopic and macroscopic states of systems, bringing entropy as a core principle in statistical physics. The transformative shift in this field began with the introduction of Shannon entropy. Shannon's \citep{shannon2001mathematical} pioneering work applied entropy to quantify the capacity and efficiency of communication channels. Since then, researchers have developed entropy measures rooted in various theoretical frameworks, such as reliability theory, temporal theory, fuzzy theory, graph theory, and fractional calculus which are applied across different fields such as artificial intelligence, decision-making, biological systems, chemical processes, communication systems, stock market, dynamical systems and signal and image processing. This makes it a versatile application-specific technique for quantifying uncertainty and disorder. Recently, Hopfield\citep{Hopfield2024} and Hinton\citep{Hinton} received the $2024$ Nobel Prize in Physics for their groundbreaking work in machine learning. Hopfield's neural memory networks and Hinton's stochastic neural network are deeply connected to entropy maximization principle to obtain system states. This recognition not only emphasizes the significance of entropy in modern computational frameworks but also motivates researchers to explore novel ideas and further applications of entropy. It underscores the need for a comprehensive review, offering researchers and practitioners a valuable resource for selecting an appropriate entropy measure to address challenges across diverse fields, ensuring both theoretical and practical relevance.

\par Several reviews of entropy measures have been conducted over the years including Shaw et al.\citep{shaw1983entropyL1}'s exploration of entropy applications in biology, economics, information science, the arts, and even religion. Golan\citep{golan2008informationL2} provided a synthesis of information theory within econometric methods, while Miracle et al.\citep{miracle2017criticalhighentropyalloythermodynamics1} highlighted critical findings in the field of high-entropy alloys. Popovic's\citep{popovic2017researchersL4} work focused on the physical interpretation of entropy, enriching the understanding of its foundational principles. Li et al.\citep{li2018entropyfaultdiagnosis2AI1} reviewed entropy-based algorithms used in fault detection. Namdari et al.\citep{namdari2019reviewL6} conducted an extensive survey on uncertainty quantification, focusing on stochastic processes and their entropy based applications. Ribeiro et al.\citep{ribeiro2021entropyL7} offered a timeline-based analysis tracing the historical development of entropy and its varied connections in different domains. However, these studies lack a detailed mathematical foundation of different entropy measures, broad application-oriented discussions, and the inclusion of essential dataset sources.

\par To address these concerns, this review traces the evolution of entropy measures, from their origins in physics to their current prominence and the applicability in diverse interdisciplinary fields. We begin by categorizing entropy measures based on their foundational theories, starting with Shannon entropy which effectively quantifies uncertainty in a system's probability distribution but fails to capture the dynamics of systems\citep{beck2009generalisedintro3} with long-range dependency, time evolution, graphical structures, non-linear systems, highly correlated processes and incomplete distributional information. The parametric generalization of Shannon entropy\citep{renyi1961measures,tsallisentropy1988possible,kaniadakis2002statisticalkaniadakisentropy} improves flexibility in modelling phenomena by deriving generalized probability distributions via principle of maximum entropy (PME), managing extremes and outliers in contingency tables, and establishing connections with other entropy measures such as Rényi and Tsallis entropy. Entropy, as a function of the probability mass/density function(pmf/pdf), captures the randomness associated with each event. In reliability analysis, the cumulative distribution function(cdf), closely linked to the survival function, offers valuable risk management and system maintenance insights\citep{rao2004cumulativeresidualentropy,di2009cumulativeEntropies,ebrahimi1996measureresidualentropy1,asadi2007dynamiccumulativeresidualentropy}. Entropy functions formulated using the cdf or survival functions establish a connection between reliability and information theory. In time-dependent systems, the entropy rate $H(X_1, X_2, \dots, X_n)/n$ analyzes the randomness of a sample $\{X_1, X_2, \dots, X_n\}$ based on their joint probability distribution, which sometimes become difficult to estimate. Capturing trends, detecting anomalies, and revealing dependency structures in time sequences have driven the development of time-sequence entropy measures such as approximate, sample, permutation, and multiscale entropy\citep{pincus1991approximateentropy,richman2004sampleentropymotivation,costa2002multiscaleentropydefination1,bandt2002permutationentropy,rostaghi2016dispersionentropy}. These are distribution-free but data-dependent statistics, making them highly versatile for measuring uncertainty without strict constraints. Human reasoning and decision-making can sometimes be subjective and incomplete, with heterogeneous and ambiguous inputs that are well-suited to modelling with fuzzy theory\citep{zimmermann2011fuzzysettheorybook}. This gives rise to fuzzy entropy\citep{zadeh1968probabilityfuzzyentropyzadeh,de1993definitionmembershipfuzzyentropy,ebanks1983measuresbrucefuzzyentropy,koskoentropyfuzzy}, a quantitative measure that provides an ordering based on vagueness and incompleteness. Extending Shannon entropy, derived from the pdf of a fuzzy set, fuzzy entropy offers a more generalized way of analyzing qualitative and quantitative data.  Fractional calculus\citep{oliveira2014reviewfractionalcalculus} is sometimes helpful to handle complex systems that exhibit self-similar structures, scale invariance, anomalous diffusion, non-locality, long-term memory. For these reasons, fractional-order entropy measures\citep{ASfractionalorderentropy,ubriaco2009entropiesUfractionalorderentropy,karci2016karcifractionalorderentropy} have been introduced to address such complexities that offers a robust modelling approach. Understanding the complexity order in graph and network relationships is crucial for applications in artificial intelligence, interdependencies, information propagation, and network optimization\citep{west2001graphtheoryintroduction}. This has led to the development of graph-theory based entropy measures\citep{1955rashevskygraphentropy,Truccographentropy1956ANO,bonchev1977informationBTgraphentropy,raychaudhury1984discriminationRRGRBgraphentropy}. 
The key contributions of this review are as follows:
\begin{enumerate}
    \item Motivation and definitions of various entropy measures with their fundamental properties.
    \item A categorical description of some of the applications of entropy measures under discussion.
    \item Sources of openly accessible data used in entropy applications along with the underlying entropy measures.
    \item Potential future directions.
\end{enumerate}


The structure of this paper is as follows: Section \ref{section2} presents definitions, connections, and properties of different entropy measures. Section \ref{section3} reviews some of the applications of entropy measures across various disciplines. Section \ref{section4} presents a thorough analysis of entropy's utility in data analysis and lists public repositories of datasets employed in these studies. Section \ref{section5} explores potential future trends, and Section \ref{section6} concludes the article.

\section{Entropy Variants: Characteristics and Properties}\label{section2}

The concept of entropy arises from the need to quantify the heat direction observed in thermodynamic processes. According to the first law of thermodynamics\citep{blundell2010concepts}, the change in internal energy $U$ of a system is given by 
\begin{equation}\label{thermodynamicEquation11}
   dU = {\dbar Q} \mskip6mu + \mskip6mu{\dbar W\vphantom{Q}},%
\end{equation}
where exact differential ${\dbar Q} \mskip6mu $ represents the heat supplied to the system, and $\mskip6mu{\dbar W\vphantom{Q}}$ is the work done on the system. While this law ensures energy conservation, it does not distinguish between reversible and irreversible processes. The second law of thermodynamics addresses this by explaining the natural flow of heat. It states that heat flows from a hotter body to a colder body as the system approaches equilibrium, and in isolation, the reverse process does not occur, as stated by Clausius. A practical illustration is the Carnot engine, where heat $Q_{\text{enter}}$ enters the system, $Q_{\text{exit}}$ leaves the system, and the work done is given by 
\begin{equation}
W = Q_{\text{enter}} - Q_{\text{exit}},
\end{equation}
when there is no change in internal energy. State variables (time-independent physical quantities) such as volume, pressure, and temperature are used to describe these phenomena mathematically. Analogously, the rate of heat entering the system per unit temperature, ${\dbar Q} \mskip6mu /T$, is an exact differential and $\int_{A}^{B}{{\dbar Q} \mskip6mu}/{T}$ is time-independent, hence qualifies as a state variable. This leads to the definition of entropy $S$, where $dS = \frac{{\dbar Q} \mskip6mu }{T}$, and thus 
\begin{equation}
S = \int \frac{{\dbar Q} \mskip6mu }{T}, 
\end{equation}
is the total heat entering the system per unit time. The relationship between entropy and internal energy from equation (\ref{thermodynamicEquation11}) is expressed as
\begin{equation}
dU = TdS - pdV,
\end{equation}
where $\mskip6mu{\dbar W\vphantom{Q}} = -pdV$ represents the work done on the system. Further in general, the laws of thermodynamics state that the internal energy of the universe, $U_{\text{universe}}$, remains constant, while the entropy of the universe, $S_{\text{universe}}$, can only increase.

Ludwig Boltzmann's pioneering work in the late 19th century established the foundations of modern statistical mechanics, linking microscopic particle behaviour to macroscopic thermodynamic properties\citep{chakrabarti2000boltzmann}. The measure of randomness in a system's distribution over its permissible microstates is described by Boltzmann's principle, known as \textbf{Boltzmann entropy} given by
\begin{equation}\label{boltzmannEntropy}
S^{W}_{B}=k_{B}\ln(W),   
\end{equation}
where, $k_B=1.380649\times10^{-23}$ (joules per kelvin) represents the Boltzmann constant, while $W$ denotes the total count of microscopic configurations consistent with the macroscopic state of the system. A function that exhibits the property of monotonically increasing and satisfies the additivity property for two independent systems uniquely characterizes the Boltzmann entropy function. Its definition is limited to calculating the uncertainty of a finite set of values. Boltzmann entropy measures the uncertainty in a dataset when the related system is in equilibrium or when all events in the sample set have equal probability. Gibbs generalized Boltzmann entropy, allowing it to capture the associated uncertainty by considering the significance of each microstate using an underlying probability distribution. Let $X$ be a discrete random variable with a corresponding probability distribution $P=\{p_i\}$. The \textbf{Gibbs entropy}\citep{Rosenkrantz2012} for the random variable $X$ is given by 
\begin{equation}\label{gibbsentropy}
    S^{X}_{G}=-k_{B}\sum_{i}p_i \ln p_i .
\end{equation}
It is a non-negative function and remains unchanged when the order of probabilities is altered. It exhibits additivity property for independent systems and is a concave function. The Gibbs entropy attains its maximum, equivalent to the Boltzmann entropy (\ref{boltzmannEntropy}), when the underlying distribution is uniform. Conversely, it reaches its minimum of zero in the case of a degenerate distribution\citep{gibbsmathematicalproperties}. \textbf{Hartley entropy}\citep{hartley1928transmission} quantifies the uncertainty related to a finite set of events or possibilities. It is defined as 
\begin{equation}\label{hartleyentropy}
    S_{H}^{X}=\ln(n) ,
\end{equation}
where $n(\geq1)$ is the size of the set of total possibilities/states/outcomes across all instances at any given time. The definition is inspired by the concept that information associated to an experiment remains constant when the number of possibilities/outcomes is the same across sample sets of different cardinality. Some of its mathematical properties include non-negativity, monotonicity with respect to the sample set size, additivity for independent systems, and achieving a minimum value when $n=1$ for a deterministic process.

A random variable $\mathcal{X}$ and its associated pdf $p_{\mathcal{X}}$ represent the system's observable and state, respectively. In a quantum system, the observable is described by a Hermitian operator $\mathcal{T}$, and the state is given by a density matrix $\mathcal{M}$. The \textbf{quantum entropy}\citep{ohya2004quantumentropy} of the state is defined as
\begin{equation}
H_{QE}(\mathcal{M}) = -\text{Tr}\left(\mathcal{M}\log\left(\mathcal{M}\right)\right),
\end{equation}
where $\text{Tr}$ computes the trace of the matrix. For any density matrix \(\mathcal{M} \), the quantum entropy \(H_{QE}(\mathcal{M}) \) is non-negative, symmetric, and concave. The minimal value is achieved when the density operator represents a pure state, which occurs when all eigenvalues possess the same eigenvector. In contrast, the highest value of \( H_{QE}(\mathcal{M}) \) is \( \log(d) \), with \( d \) representing the system's dimension, attained in the maximal mixed state.

In 1948, Shannon made a groundbreaking contribution to information theory and communication engineering. His work not only established the field of information theory but also set the foundational principles for the efficient and dependable transmission of information in modern communication systems, setting the stage for the digital era. If  $p_1, p_2, ..., p_n$ denotes the probabilities associated with the possible occurrences, then the \textbf{Shannon entropy}\citep{shannon2001mathematical} $H$ is given by
\begin{equation} \label{shannonentropy}
    H(p_1, p_2, ..., p_n)=-\sum_{i=1}^{n}p_i \ln p_i.
\end{equation}
It is clear from equation (\ref{gibbsentropy}) that Shannon entropy varies in proportion to Gibbs entropy. The uniqueness of the functional form in equation (\ref{shannonentropy}) is established by satisfying sufficient conditions: $H$ is continuous with respect to each $p_i$'s, monotonically increases with the sample size under the uniform distribution, and maintains a weighted sum of individual $H$ values when choices are decomposed into successive steps. For a continuous random variable $X$ with pdf $p(x)$, the Shannon entropy\citep{continuousshannonentropycover1999elements} in continuous version is defined as
\begin{equation}\label{continuousShannonEntropy}
H_{X}=-\int_{-\infty}^{\infty}p(x)\ln(p(x))dx.
\end{equation}

\begin{figure}
    \centering
    \includegraphics[width=1\linewidth]{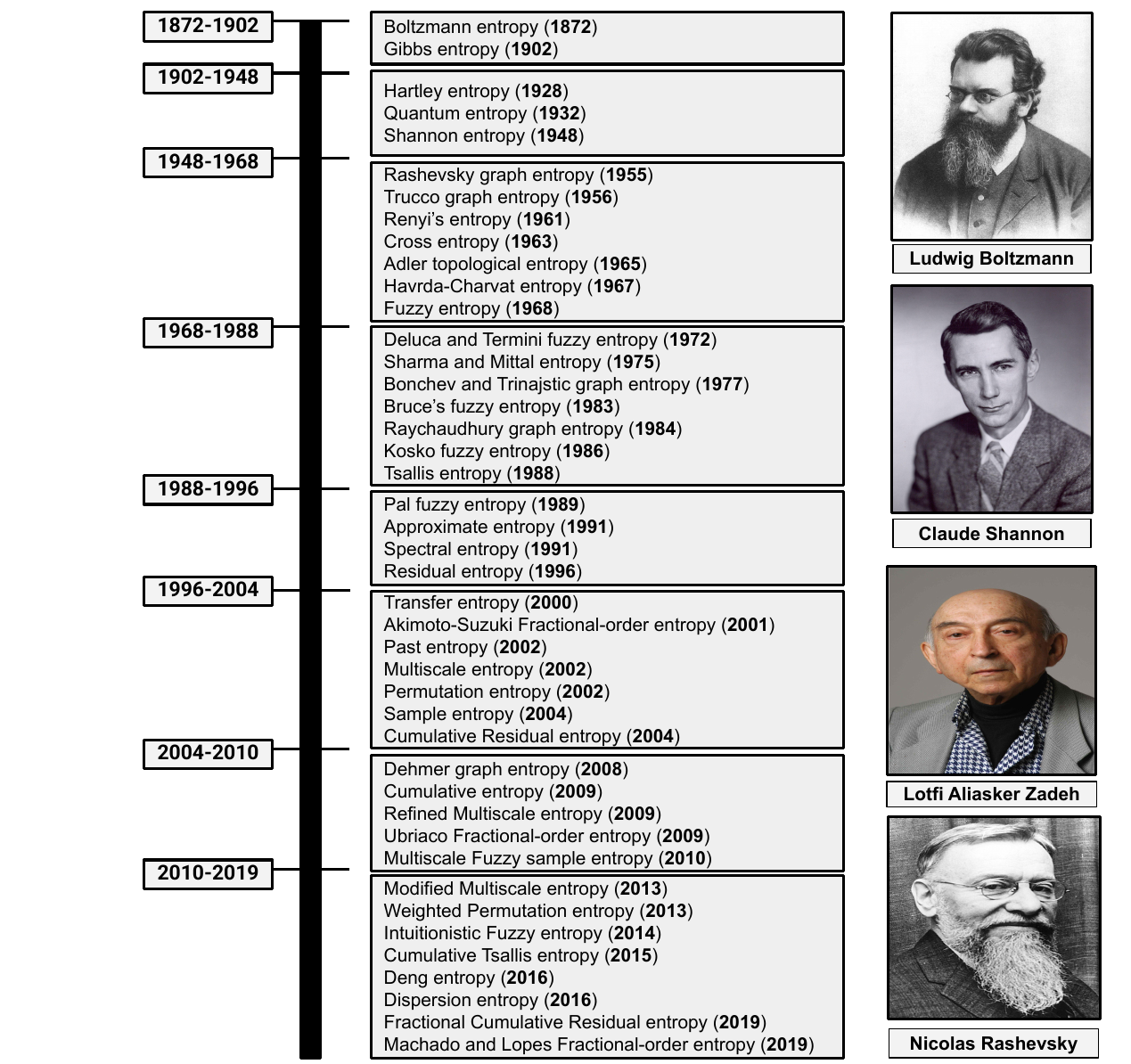}
    \caption{Evolution of entropy measures over years}
    \label{fig:8}
\end{figure}

This formula is the continuous analogue of the Shannon entropy defined in equation (\ref{shannonentropy}), commonly referred to as \textbf{Differential entropy}. One significant distinction between the discrete and continuous versions of Shannon entropy is that for the discrete random variable, Shannon entropy is always positive, but for a continuous random variable, as given by equation (\ref{continuousShannonEntropy}), it can also be a negative real number. The commonly utilized mathematical properties of both measures are additivity, expansibility, concavity, Lesche-stability, composability, and attained the maximum value under uniform distribution. Khinchin\citep{khinchinUniqueness1957} provides four axioms that establish the uniqueness of the functional form of Shannon entropy (\ref{shannonentropy}). These axioms are as follows: 
\begin{itemize}
    \item[1.] $H(\{p_i\})$ is continuous in $p_i$'s.
    \item[2.] For uniform distribution, $H$ increases as the sample size increases.
    \item[3.] $H(p_1,p_2,...,p_n,0)$ and $H(p_1,p_2,...,p_n)$ are equal.
    \item[4.] $H(X,Y)=H(X|Y)+H(Y)$, for the joint random variables $(X,Y)$.
\end{itemize}
For the flexibility and thorough analysis of complex systems, numerous generalizations of Shannon entropy exist, including their continuous versions. Next, we summarize several parametric generalizations of the widely used Shannon entropy measure, which involve modifications to the axioms mentioned in (1)-(4).

\subsection{Parametric generalizations of Shannon Entropy}
A one-parametric entropy measure consistent with the additivity property for generalized distributions is defined by Rényi\citep{renyi1961measures} using the nonlinear average. In the case of a discrete random variable $X$ with distribution \( P \), \textbf{Rényi's entropy} of order $\alpha>0(\neq 1)$ is defined as 

\begin{equation}\label{renyientropydiscrete}
    H_{X}^{\alpha}=\frac{1}{1-\alpha}\ln \left( \sum_{i=1}^{n} p_i^{\alpha} \right),
\end{equation} 
and for a continuous random variable $X$ with pdf \( p(x) \), it is defined as
\begin{equation}\label{renyientropycontinuous}
    H_{X}^{\alpha}=\frac{1}{1-\alpha}\ln \Biggl( \int_{-\infty}^{\infty} p^{\alpha}(x)dx\Biggr).
\end{equation} 
It can be seen that (\ref{shannonentropy}) and (\ref{continuousShannonEntropy}) are limiting cases of (\ref{renyientropydiscrete}) and (\ref{renyientropycontinuous}) as \( \alpha \to 1 \). Rényi's entropy for a discrete distribution is always a non-negative real number, whereas, for a continuous distribution, it can also be negative\citep{tabass2016renyientropyIs_not_always_non_negative}. Additionally, the continuous form (\ref{renyientropycontinuous}) of Rényi's entropy is not the limit of the discrete form (\ref{renyientropydiscrete}), in which the probabilities of discrete random variable estimated from the creation of bins of appropriate lengths. Jizba and Arimitsu\cite{renyientropy_properties} derived the functional uniqueness theorem for (\ref{renyientropydiscrete}) based on a set of five axioms. These axioms are: continuity of a function $H_{X}^{\alpha}$ with respect to each $p_i$'s, maximum at uniform distribution, expansibility, $H_{X,Y}^{\alpha} = H_{Y|X}^{\alpha} + H_{X}^{\alpha}$, where 
\begin{equation}
H_{Y|X}^{\alpha}=f^{-1}\left( \sum_{k=1}^{n}g_{k}(\alpha)f\left( H_{Y|X=x_k}^{\alpha} \right) \right)
\end{equation}
with $g_{k}(\alpha)={{(p_k)}^{\alpha}}/\sum_{k=1}^{n}{{(p_k)}^{\alpha}}$ is generalized average, and $f$ is invertible and positive on $[0, \infty)$. 

Havrda and Charvat proposed a quantitative measure for the classification problem of a non-empty set $B$ with a normed measure $\mu$, also known as \textbf{Havrda–Charvat entropy}. Here, the measure $\mu$ need not be a probability measure. Let $M_1, M_2, \cdots, M_n$ be a partition of the non-empty set $B$ with normed measure values $\mu_i$ for each $M_i$. Then, the Havrda–Charvat entropy\citep{havrda1967quantificationHarvatCharvatentropy} is defined as
\begin{equation}\label{Havrda_entropy}
    S_{B}^{a}=\frac{2^{a-1}}{2^{a-1}-1} \left( 1-\sum_{i=1}^{n}\mu_i^{a} \right),
\end{equation}
where $a(\neq 1)>0$. As the parameter $a$ approaches $1$, the function form in (\ref{Havrda_entropy}) tends to (\ref{shannonentropy}), which corresponds to Shannon entropy. This makes $S_{B}^{a}$ one parametric generalization of Shannon entropy. The construction of $S_{B}^{a}$ is motivated by four axioms that characterize it as a unique function: continuity of $S_{B}^{a}$ for each $\mu_i$ such that $\sum_{i=1}^{n}\mu_i=1$ for all $a(\neq 1)>0$, $S_{B}^{a}$ is zero for degenerate distributions, equals $1$ for the distribution $\{ 1/2, 1/2 \}$, $S_{B}^{a}$ is expansible, and 
\begin{equation}
    S_{B}^{a}(\mu_1,\cdots,\mu_{i-1},\nu_{i_1},\nu_{i_2},\mu_{i+1},\cdots,\mu_{n})=S_{B}^{a}(\mu_1,\cdots,\mu_{n})+{\beta}{\mu_i^a} S_{B}^{a}\left(\frac{\nu_{i_1}}{\mu_i},\frac{\nu_{i_2}}{\mu_i}\right),
\end{equation}
for every $i=1,2,\cdots, n$, $\nu_{i_1}+\nu_{i_2}=\mu_i >0$ and $\beta>0$. In the region of $\mu_i \geq 0$ for all $i=1,2, \cdots, n$ and $\sum_{i=1}^{n}\mu_i=1$, $S_{B}^{a}$ is concave and achieves its maximum value at the uniform distribution. For a probability distribution $\{p_i \}$, it is easy to see the relation 
\begin{equation}\label{relation between renyi's and Havrda entropy}
    S_{B}^{a}(p_1, \cdots, p_n)=\frac{2^{a-1}}{2^{a-1}-1}\Bigg(1- \exp\left((1-a)H_{X}^{a}(p_1, \cdots, p_n)\right) \Bigg),
\end{equation}
between the Havrda-Charvat (equation (\ref{Havrda_entropy})) and discrete Renyi's entropy (as in equation (\ref{renyientropydiscrete})), here the random variable $X$ is defined on the partition of non-empty set $B$. 

Sharma and Mittal\citep{sharmaandmittalentropy1975new} proposed a two-parameter extension of (\ref{shannonentropy}). This new approach generalizes the additivity property and allows greater flexibility through the generalized average due to the utility of the escort distribution. Let $X$ be a discrete random variable then \textbf{Sharma and Mittal entropy} is given as
\begin{equation} \label{sharmaandmittalentropy_discrete}
    H_{\alpha,\beta}^{X}=\frac{1}{1-q}\left( {\left(\sum_{i}^{n}p_i^{\alpha} \right)}^{\frac{1-\beta}{1-\alpha}} -1\right),
\end{equation}
where $\alpha(\neq 1)$ is non-negative and $\beta(\neq 1)$ is any real number. It is easy to see that if $\beta$ approaches $1$ with $\alpha\neq 1$, then $H_{\alpha,\beta}^{X}$ converges to (\ref{renyientropydiscrete}), and if both $\alpha$ and $\beta$ simultaneously approach $1$, then $H_{\alpha,\beta}^{X}$ converges to (\ref{shannonentropy}). Here are some important properties:
\begin{enumerate}
    \item $H_{\alpha,\beta}^{X}$ is maximum at uniform distribution.
    \item $H_{\alpha,\beta}^{X}$ is a continuous function for each $p_i$'s.
    \item $H_{\alpha,\beta}^{X}$ is expansible.
    \item $H_{\alpha,\beta}^{X}(1/2,1/2)=\log_{q}(1)$, where 
          \begin{equation}
              \log_{q}(x)=\begin{cases}
                             \log(x), & \text{if $q=1$,} \\
                             \frac{x^{1-q}-1}{(1-q)\ln(2)}, & \text{ if $q\neq 1$.}
                         \end{cases}
          \end{equation}
    \item Let $R=\{r_{ij} \}_{i=1,j=1}^{n,m}$ be a joint pmf of $(X, Y)$. Also, let $p_i=\sum_{j=1}^{m}r_{ij}$, $q_j=\sum_{i=1}^{n}r_{ij}$, and $\alpha(\geq0)$ be a real constant. If the distribution of $(Y|X=x_k)$ is $Q_{|k}=\{q_{j|k}\}_{j=1}^{m}$, where $q_{j|k}={r_{kj}/p_k}$, then,
    \begin{equation}
        H_{\alpha,\beta}^{X,Y}=H_{\alpha,\beta}^{X}{\oplus}_{k} H_{\alpha,\beta}^{Y|X},
    \end{equation}
    where 
    \begin{equation}
        H_{\alpha,\beta}^{Y|X}=g^{-1}\left( \sum_{i=1}^{n}p_{i}^{(\alpha)}g\left(H_{\alpha,\beta}^{Y|X=x_k} \right) \right),
    \end{equation}
    for any real numbers $a_1,b_1,k_1$,
    \begin{equation}
        a_1{\oplus}_{k_1}b_1=a_1+b_1+(1-k_1)a_1b_1,
    \end{equation}
    and $g$ is a continuous invertible function, and  $p_{i}^{(\alpha)}=\frac{p_{i}^{\alpha}}{\sum_{k=1}^{n}p_{k}^{\alpha}}$ is the escort distribution.
\end{enumerate}
Note that if a function satisfies all the properties from $1-5$, it can be uniquely defined\citep{ilic2021alphaSharmaAndMittalproperty}, as in equation (\ref{sharmaandmittalentropy_discrete}). For a continuous random variable $X$ with pdf $p(x)$, it can be defined as 
\begin{equation}
H_{\alpha,\beta}^{X}=\frac{1}{1-\beta}\left({\left(\int_{-\infty}^{\infty}{p^{\alpha}(x)} \right)}^{\frac{1-\beta}{1-\alpha}} -1\right).
\end{equation}

Tsallis\citep{tsallisentropy1988possible} generalizes Shannon entropy into a one-parameter form based on the multifractal scaling of a quantity. The \textbf{Tsallis entropy} for a discrete distribution $P$ and a real parameter $q(\neq 1)$ is defined as 
\begin{equation} \label{Tsallisentropydiscretecase}
    T_q^{X}=\frac{k}{q-1}\left(1-\sum_{i=1}^{n}p_{i}^{q} \right),
\end{equation}
here $k$ is a conventional positive constant. The $T_q^X$ function is always non-negative for a discrete random variable $X$ and entropic index $q(\neq 1)$. Also, it can take negative values for continuous distributions, depending on the value of $q$, such as in the case of the exponential distribution\citep{kumar2024estimation}. It reaches its maximum (minimum) at a uniform distribution when $q>0(q<0)$. Additionally, the function has the expansibility property, and $T_q^X$ is concave (convex) for $q > 0$ ($q < 0$). One notable property\citep{tsallisentropyproperties2009introduction} distinguishing Tsallis entropy from other entropy functions is its non-additivity. Specifically, for the joint distribution of independent random variables $X$ and $Y$, the Tsallis entropy for $k=1$ is given by 
\begin{equation}\label{nonadditivityofTsallisentropy}
    T_q^{(X,Y)}=T_q^{X}+T_q^{Y}+(1-q)T_q^{X}T_q^{Y}.
\end{equation}
The cases \( q < 1 \) and \( q > 1 \) are often defined as superadditive and subadditive, respectively. It also has a connection with the Jackson derivative and satisfies the definition of experimental robustness for \( q > 0 \). Other notable variations\citep{tsallisentropyproperties2009introduction} of Tsallis entropy are the escort Tsallis entropy, defined as 
\begin{equation}\label{EscortTsallisentropy}
    T_q^{E}(X)=\frac{k}{q-1}\left(1-\sum_{i=1}^{n}\frac{p_{i}^{q}}{\sum_{j=1}^{n}p_j^q} \right),
\end{equation}
and the normalized Tsallis entropy, given by 
\begin{equation}\label{Normalizedtsallisentropy}
    T_q^{N}(X)=\frac{T_q^{X}}{\sum_{i=1}^{n}p_i^{q}}.
\end{equation}
The continuous version of Tsallis entropy is defined as
\begin{equation}\label{continuoustsallisentropy}
    T_q^{X}=\frac{k}{q-1}\left(1-\int f_X^q(x)dx \right),
\end{equation}
where $q(\neq1)>0$. Similar to Shannon entropy for continuous cases, this can also exhibit negative values. Shannon entropy is a particular case of Tsallis entropy; as $q$ goes to $1$, Tsallis entropy tends to Shannon entropy. Additionally, the well-known relationship between Tsallis (\ref{Tsallisentropydiscretecase}) and Rényi entropy (\ref{renyientropydiscrete}) for $q(\neq 1)>0$, is given by 
\begin{equation} \label{relationbetweenRenyiandTsallisentropy}
    H_{X}^{q}=\left(\frac{1}{1-q}\right)\ln\left(1-(q-1){T_{q}^{X}} \right), \ \ \ \text{ for $k=1$.}
\end{equation}
From equation (\ref{sharmaandmittalentropy_discrete}), it is evident that as \( \alpha \) approaches \( \beta \), the Sharma and Mittal entropy converges to the Tsallis entropy (\ref{Tsallisentropydiscretecase}). The relationship between the Havrda–Charvat entropy (\ref{Havrda_entropy}) and Tsallis entropy (\ref{Tsallisentropydiscretecase}) is expressed as 
\begin{equation}
    S_B^a=\left(\frac{a-1}{1-2^{1-a}}\right)T_X^a, \text{\ \ \ where $a(\neq1)>0$.}
\end{equation}
\textbf{Kaniadakis entropy}, another generalization of Shannon entropy, emerged to describe deformations in Einstein's special relativity theory\citep{kaniadakis2002statisticalkaniadakisentropy}. It is defined as
\begin{equation}\label{kaniadakis entropy}
    H_{X}^{k}=-\sum_{i}\frac{p_i^{1+k}-p_i^{1-k}}{2k},
\end{equation}
where $k$ is a non-zero real parameter. If \( k \) approaches 0, \( H_X^k \) tends to \( H \) (\ref{shannonentropy}), Shannon entropy. In the definition of Kaniadakis entropy, a one-parametric generalization of the natural logarithm is used, given as 
\begin{equation}\label{k-log(kaniadakisentropy)}
    \ln_k(x)=\frac{x^k-x^{-k}}{2k},
\end{equation}
which converges to the natural logarithm as $k\to0$. The $k$-logarithm function is concave for $k\in(-1,1)$, which ensures the concavity of Kaniadakis entropy. Additionally, it increases with \( x \) monotonically.  Let 
\begin{equation}
\mathcal{u}_k(x)=\frac{x^k-x^{-k}}{2k},
\end{equation}
and 
\begin{equation}
    I_{k}(X)=E(\mathcal{u}_k(X)).
\end{equation}
The additivity property for Kaniadakis entropy is given by
\begin{equation}
    H_{X,Y}^{k}=I_{k}(X)H_{Y}^{k}+I_{k}(Y)H_{X}^{k},
\end{equation}
where the random variable $X$ is independent of $Y$. Moreover, the inverse of the \( k \)-logarithm, given in (\ref{k-log(kaniadakisentropy)}), is the \( k \)-exponential, defined by 
\begin{equation}
    \exp_{k}(x)={\left(\sqrt{1+{kx}^2}+kx \right)}^{1/k}.
\end{equation}
This is a positive valued convex function that monotonically increases and converges to the standard exponential function as \( k \) approaches $0$.

These are some widely accepted parametric generalizations of Shannon entropy. Next, we explore diverse prominent entropy functions formulated from the cdf.

\subsection{Entropy Functions Derived from the Cumulative Distribution Function}

\begin{sidewaystable} 
\centering
\begin{tabular}{p{1.5cm}p{2cm}p{4.5cm}p{8.5cm}p{7.5cm}}
\toprule
\textbf{S. No.} & \textbf{Reference} &\textbf{Name of entropy}&\textbf{Defination of entropy} &\textbf{Distinctive properties} \\
\midrule
 $1$ & \citep{rajesh2019cumulativetsallisentropy} & Cumulative Tsallis entropy & 
 $H_{CTE}^{\alpha}(X)=\frac{1}{\alpha-1}\int_{0}^{\infty}\left(\overline{F}_X(x)-{(\overline{F}_X(x))}^{\alpha} \right)dx$, \ \ \  $0<\alpha\neq 1$
 & $\lim_{\alpha\to1}H_{CTE}^{\alpha}(X)=E_{CRE}(X)$ and for $0<\alpha<1( \alpha>1)$, $H_{CTE}^{\alpha}(X)\geq(\leq)0$.  \\
 $2$ & \citep{mirali2017weightedcumulativeresidualentropy} & Weighted cumulative residual entropy & $H_{WCRE}(X)=-\int_{0}^{\infty}x\overline{F}_X(x)\log(\overline{F}_X(x))dx$ &  If $p>2$ such that $E(X^p)<\infty$ then $H_{WCRE}(X)<\infty$. $H_{WCRE}(X)=0$ iff $X$ is degenerate.\\
 $3$ & \citep{xiong2019fractionalcumulativeresidualentropy} & Fractional cumulative residual entropy & $H_{FCRE}^{q}(X)=\int_{0}^{\infty}\overline{F}_X(x)\left[-\log(\overline{F}_X(x))\right]^{q}dx$, $0\leq q \leq 1$ & $H_{FCRE}^{1}(X)=E_{CRE}(X)$, $H_{FCRE}^{q}(X)$ is non-negative and non-additive, it is convex with respect to $q$ and concave function of cdf. \\
 $4$ & \citep{sunoj2012dynamiccumulativeresidualrenyisentropy} & Dynamic cumulative residual Renyi’s entropy & $H_{DCRRE}^{\alpha}(X;t)=\frac{1}{1-\alpha}\log\left( \int_{t}^{\infty}\frac{\overline{F}^{\alpha}_X(x)}{\overline{F}^{\alpha}_X(t)}dx\right)$, $0<\alpha \neq 1$ & $H_{DCRRE}^{\alpha}(X;t)$ uniquely determines $\overline{F}^{\alpha}_X(x)$ and for a absolutely continuous survival function $S$ and hazard rate $h_{X}(t)$, we have $(1-\alpha){H_{DCRRE}^{\alpha}}^{'}(X;t)=ch_{X}(t)$, where $c$ is a constant. \\
 $5$ & \citep{sati2015somedynamiccumulativeresidualtsallisentropy} & Dynamic cumulative residual Tsallis entropy & $H_{DCRTE}^{q}(X;t)=\frac{1}{q-1}\left(1-\int_{t}^{\infty}{\left(\frac{\overline{F}^{\alpha}_X(x)}{\overline{F}^{\alpha}_X(t)} \right)}dx \right)$, $0<q\neq 1$ & If \( h_X(t) \leq h_Y(t) \), then \( H_{DCRTE}^{q}(X;t) \leq (\geq) H_{DCRTE}^{q}(Y;t) \) for \( q>1 (0<q<1)\). \\
 $6$ & \citep{di2002pastentropy} & Past entropy & $H_{PE}(X;t)=-\int_{0}^{t}\frac{f_{X}(x)}{F_{X}(t)}\log \frac{f_{X}(x)}{F_{X}(t)}dx$ & Shannon entropy $H_{X}=F(t)H_{PE}(X;t)+\overline{F}_X(t)H_{Res}(f_{T};T)+c_1$, where $c_1$ is the shannon entropy of  Bernoulli distribution. \\

\bottomrule
\end{tabular}
\caption{Functions inspired from the cumulative and residual entropy}
\label{Table of cumulative entropies}
\end{sidewaystable}

The parametric generalization of Shannon entropy relies solely on the distribution of discrete (or continuous) random variables through their pmf or pdf. However, deriving the pmf or pdf of a random variable, for example, when it involves a mixture of Gaussian and delta distributions, can be challenging. Further, the other irregularity is that while Shannon entropy is always positive for a discrete distribution, it can assume negative values for continuous random variables. 

\noindent The \textbf{Cumulative residual entropy}(CRE)\citep{rao2004cumulativeresidualentropy} derives from the cdf of random variables. Let \( \textbf{X}\in \mathbb{R}^n \) be a random vector then CRE is given by 
\begin{equation}\label{cumulativeresidualentropy}
    E_{CRE}(\textbf{X})=-\int_{\mathbb{R}^n_{+}}P\left(|\textbf{X}|>x\right)\log P\left(|\textbf{X}|>x\right)dx
\end{equation}
where \( \textbf{X} = (X_1, X_2, \cdots, X_n) \), \( x = (x_1, x_2, \cdots, x_n) \), and \( |\textbf{X}| > x \) implies \( |X_i| > x_i \). The CRE is a non-negative, concave function that can also be computed from sample data easily, converging to the exact value asymptotically. Also, if there exists \( n_1 > n \) such that \( E(|X_i|^{n_1}) < \infty \) for all \( i \), then \( E_{CRE}(\textbf{X}) < \infty \). For non-negative independent random vectors \( \textbf{X} \) and \( \textbf{Y} \), \( \max\{E_{CRE}(\textbf{X}), E_{CRE}(\textbf{Y})\} \leq E_{CRE}(\textbf{X}+\textbf{Y}) \). The interrelation between the CRE and Shannon entropy (\ref{continuousShannonEntropy}), for a random variable \( X \), characterized by its pdf \( f \) is expressed as 
\begin{equation}
    E_{CRE}(X)\geq {c_0}H_X,
\end{equation}
where $c_0 	\cong 0.2065 $.

A complementary measure, functioning as a dual to CRE and based on cdf, is \textbf{Cumulative entropy}(CE). The cumulative entropy\citep{di2009cumulativeEntropies} for a non-negative random variable \( X \) is defined as 
\begin{equation} \label{cumulativeentropy}
    CE(X)=-\int_{0}^{\infty}F(x)\log F(x) dx,
\end{equation}
where, \( F \) is the cdf of $X$. CE is non-negative and equals zero for a degenerate random variable. Note that, for a symmetric cdf \( F \), with respect to the finite mean of \( X \), the CE equals the CRE. Following are some important results related to CE:
\begin{enumerate}
    \item If $a>0$ and $b\geq 0$ then we have $CE(aX+b)=aCE(X)$.
    \item If \( a \) and \( b \) are finite real numbers such that the random variable \( X \) with support \([0, a]\) is independent of \( Y \) with supports \([0, b]\) respectively, then we have \(CE(X,Y) = [b - E(Y)]CE(X) + [a - E(X)]CE(Y) \).
    \item The relation 
    \begin{equation}
        CE(X)=E[G(X)],
    \end{equation}
    where $G(x)=-\int_{x}^{\infty}\log F(y)dy$, $x\geq 0$ 
holds for a finite CE of an absolutely continuous random variable.
\item If the random variable $X$ is absolutely continuous then we have
\begin{equation}
    CE(X)\geq c_{0}\exp\left(H_{X}\right).
\end{equation}
\item When \( X \) and \( Y \) are independent random variables with the support of non-negative real numbers, the inequality \( \max \{CE(X), CE(Y) \}\leq CE(X+Y) \) holds.
\end{enumerate}
The information and uncertainty inherent in a lifetime distribution are studied using residual entropy. The \textbf{Residual entropy}\citep{ebrahimi1996measureresidualentropy1} of a non-negative random variable \(\mathcal{T}\), given the survival until time \(t\) is defined as
\begin{equation}\label{residualentropy}
\begin{split}
    H_{Res}(f_{\mathcal{T}};t)&=-\int_{t}^{\infty}\frac{f_{\mathcal{T}}(x)}{\overline{F}_{\mathcal{T}}(t)}\log\left(\frac{f_{\mathcal{T}}(x)}{\overline{F}_{\mathcal{T}}(t)}\right)dx \\
    & = 1-\frac{1}{\overline{F}_{\mathcal{T}}(t)}\int_{t}^{\infty}\log\left( \mathcal{L}_{F_{\mathcal{T}}}(x) \right)f_{\mathcal{T}}(x)dx,
\end{split}
\end{equation}
where $F_{\mathcal{T}}$ is the cdf of $\mathcal{T}$, $\overline{F}_{\mathcal{T}}=1-F_{\mathcal{T}}$ and $\mathcal{L}_{F_{\mathcal{T}}}(t)=f_{\mathcal{T}}(t)/F_{\mathcal{T}}(t)$. It can take negative values. Also, at \(t = 0\), 
\begin{equation}
    H_{Res}(f_{\mathcal{T}};0)=H_{\mathcal{T}},
\end{equation}
which is continuous Shannon entropy as in equation (\ref{continuousShannonEntropy}). Note that if $H_{Res}(f_{\mathcal{T}};t)$ is finite for \(t \geq 0\) and $f_{\mathcal{T}}(t)$ is a continuous function, then $H_{Res}(f_{\mathcal{T}};t)$ uniquely determines $\overline{F}_{\mathcal{T}}(t)$. Another result states that if \(d_{F}(t)\) represents the mean time until failure, given survival up to time \(t\), calculated as \(E(T - t \mid T > t)\), then for a finite \(d_{F}(t)\), it holds that \( H_{Res}(f_{\mathcal{T}};t)\leq 1 + \log(d_{F}(t))\).
The entropy derived from the residual distribution (given in (\ref{residual_distribution})) with the survival function is known as \textbf{Dynamic cumulative residual entropy} to explore the uncertainty and information associated with the residual lifetime. Let \(\overline{F}(x)\) be the residual distribution associated with the non-negative random variable \(X\), and define the function 
\begin{equation} \label{residual_distribution}
\bar{F}(x;t) = \begin{cases}
    \frac{\overline{F}(x)}{\overline{F}(t)}, & \text{if $x>t$,} \\
    $1$, & \text{otherwise,}\
\end{cases} 
\end{equation}
where, $t(\geq0)$ is a real number. Then, the dynamic cumulative residual entropy\citep{asadi2007dynamiccumulativeresidualentropy} is defined as
\begin{equation}\label{dynamiccumulativeresidulentropy}
    H_{DCRE}(X;t)=-\int_{t}^{\infty}\overline{F}(x;t)\log(\overline{F}(x;t))dx.
\end{equation}

For \(t = 0\), the dynamic cumulative residual entropy equals the CRE (\ref{cumulativeresidualentropy}). According to the definition, the cdf \(F\) is termed as decreasing (increasing) dynamic cumulative residual entropy if \(H_{DCRE}(X;t)\) is a decreasing (increasing) as a function of \(t\). This gives, \(F\) is decreasing (increasing) dynamic cumulative residual entropy if and only if \(H_{DCRE}(X;t) \leq (\geq) d_{F}(t)\) for \(t \geq 0\). $E({(X-t)}^2|X>t)/2d_{F}(t)$ is the maximum value that $H_{DCRE}(X;t)$ can attain, and this upper limit is attained if and only if $X$ is an exponentially distributed random variable. We also present some widely used entropy measures inspired by cumulative and residual entropies in Table \ref{Table of cumulative entropies}.

\subsection{Temporal entropy Measures}
Temporal-based entropy measures are crucial in understanding the dynamics of time-dependent systems. Unlike traditional entropy measures, which assess the disorder or randomness in a static dataset, temporal-based entropy considers the evolution of data over time, capturing the complexity and unpredictability of temporal sequences. It helps in identifying patterns, trends, and irregularities within the time-dependent systems.

Regularity in a time-dependent system refers to the consistency and predictability of patterns within data. \textbf{Approximate Entropy} (AE) is a measure of regularity in temporal data. Given a sequence of data \(\{ a_1, a_2, \ldots, a_N \}\), define \( b_i = \{ a_i, \ldots, a_{i+m-1} \} \) where \( r(>0) \) be a real constant and \( m \) be a non-negative integer less than \( N \). The distance between the tuples \( b_i \) and \( b_j \) is calculated as 
\begin{equation}
d(b_i, b_j) = \max_{k=1,...,m}\left(|a_{i+k-1}-a_{j+k-1}|\right).
\end{equation}
Let
\begin{equation}
    p_{m,r}(i)=\frac{\#\left(j;j\leq N-m+1 \text{ \ \ and \ \ } d(b_i, b_j)\leq r \right)}{(N-m+1)}
\end{equation}
and 
\begin{equation}
    \pi(m,r)=\frac{1}{N-m+1}\sum_{i=1}^{N-m+1}\log(p_{m,r}(i)).
\end{equation}
Then, the AE\cite{pincus1991approximateentropy} estimator is 
\begin{equation}
    H_{AE}(m,r;N)=\pi(m,r)-\pi(m+1,r),
\end{equation}
and defined by
\begin{equation}
  H_{AE}(m,r)= \lim_{N\to \infty} H_{AE}(m,r;N).
\end{equation}
Recommended values\citep{delgado2019approximateentropypropertiesSamplealso} for the parameter \( m \) are typically low, such as $2$ or $3$. The parameter \( r \) should be sufficiently large to avoid trivial cases. The value of $m$ is generally small; even a dataset with as few as $N=100$ points is sufficient for analysis. AE can be used without the need for prior information or assumptions about the dataset or the underlying process of generating the values.

The two important limitations of AE are the dependency of the outcome on $r$ and the cardinality of the temporal dataset. An attempt in the form of \textbf{Sample Entropy} \citep{richman2004sampleentropymotivation} is made to overcome these limitations. For a given value of \( r \), \( p_{m,r}(i) \) represents the conditional probability of \( m \)-length tuples matching or being close to \( b_i \) among all \( (n-m+1) \) \( m \)-length tuples (\( b_j \)) of data points, based on the \( r \) value, with self-matching ensuring non-emptiness. If \( S_i^1 \) is the set of all possible \( m \)-length tuples and \( S_i^2 \) is the set of matching \( m \)-length tuples, then AE calculates the ratio \((|S_i^2| + 1) / (|S_i^1| + 1)\). It produces a bias since the correct ratio should be \(|S_i^2| / |S_i^1|\), which affects the results, especially when dealing with small sample sizes. Sample Entropy \citep{richman2000physiologicalsampleentropydefination} is computed as follows.
Let
\begin{equation}
    p_{m}^{r}(i)=\frac{1}{N-m-1}\left(\text{Number of $b_j$ at maximum distance of $r$ to $b_i$ except itself}\right)
\end{equation}
and 
\begin{equation}
    B_{m}^{r}=\frac{1}{N-m}\sum_{i=1}^{N-m}p_{m}^{r}(i).
\end{equation}
Similarly, if we consider \( (m+1) \)-length tuples, denoted as \( b_i \), we can define the terms 
\begin{equation}
    p_{m+1}^{r}(i)=\frac{1}{N-m-2}\left(\text{Number of $b_j$ at a distance of $r$ to $b_i$ except itself}\right)
\end{equation}
and 
\begin{equation}
    A_{m}^{r}=\frac{1}{N-m-1}\sum_{i=1}^{N-m-1}p_{m+1}^{r}(i).
\end{equation}
The estimator of sample entropy is 
\begin{equation} \label{sampleentropywiththreeparameter}
    H_{SE}(m,r,N)=-\log\left( \frac{A_{m}^{r}}{B_{m}^{r}} \right)
\end{equation}
and the sample entropy (SE) is defined as 
\begin{equation}
    H_{SE}(m,r)=\lim_{N\to \infty}\left[-\log\left( \frac{A_{m}^{r}}{B_{m}^{r}} \right) \right].
\end{equation}
It is always non-negative. Generally, we consider \( r = \) (standard deviation of the time sequence)$/5$ and \( m = 2 \). A low SE value indicates regularity in temporal data, whereas a high value suggests irregularity. By varying the parameter \( m \), short-term and long-term patterns in the time sequence can be revealed. For regular and periodic data, the SE equals zero; for uncorrelated random data, the maximum entropy value is attained \citep{delgado2019approximateentropypropertiesSamplealso}.

Coarse-graining\citep{humeau2015multiscalecoarsegraineddefination} in time series analysis is valuable for simplifying complex data and, uncovering significant patterns and trends over different time scales. A coarse-graining method produces a series of time sequences representing the system's behaviour at various time scales. The average of the data points within consecutive, non-overlapping intervals of length \( k \) forms the coarse-grained time series for a given scale \( k \). Thus, given a univariate time sequence \(\{x_i\}\) of length \( M \), the coarse-grained time series \(\{y^d_j\}\) of scale \( d \) is calculated as 
\begin{equation}\label{multiscaleentropycoarsegrainedseriesdefination}
    y^d_j=\frac{1}{d}\sum_{i=(j-1)d+1}^{jd}x_i, \ \ \ 1\leq j \leq M/d.
\end{equation}
It is evident that at scale one, both coarse-grained and the original time sequence are identical. The computation of \textbf{Multiscale Entropy} \citep{costa2002multiscaleentropydefination1,costa2005multiscaleentropydefination2} involves two steps: 
\begin{enumerate}
    \item Given a scale \( d \) (where \( M/d \) is a positive integer), calculate the coarse-grained time sequence \(\{y^d_j\}\) and
    \item Calculate the SE using the above coarse-grained time sequence.
\end{enumerate}
The choices of the scale parameter depend on the domain knowledge, data length, and empirical testing while also considering the need to avoid overfitting and ensure accurate trend interpretation. At scale \( d = 1 \), SE can be considered a particular multiscale entropy case. Both entropy functions share similar mathematical properties, including function definition, non-negativity, sensitivity to dependence parameters, and noise sensitivity. The coarse-graining in equation (\ref{multiscaleentropycoarsegrainedseriesdefination}) employed here is similar to applying a finite impulse response filter\cite{humeau2015multiscalecoarsegraineddefination} given as 
\begin{equation}
    y^d_j=\frac{1}{d}\sum_{l=0}^{d-1}x_{(j-l)}, \ \ \ 1\leq j \leq M,
\end{equation}
to the time sequence and downsampling by a factor of \( d \). Table \ref{Table of generalized entropy based on multiscale entropy} presents various known generalized multiscale entropy functions, along with their definitions and distinguishing properties.


\begin{sidewaystable} 
\centering
\begin{tabular}{p{1cm}p{1.2cm}p{3cm}p{5.5cm}p{7.5cm}p{4.5cm}}
\toprule
\textbf{S. No.} & \textbf{Reference} &\textbf{Name of entropy}&\textbf{Coarse-graining Procedure}&\textbf{Defination of entropy} &\textbf{Distinctive properties} \\
\midrule
 $1$ & \citep{chen2009measuringfuzzymultiscalesampleentropy} & Multiscale Fuzzy Sample Entropy 
 &  No Coarse-graining
 &  For \( m \)-tuple vectors \( b_i \), consider vectors \( b_{i}-b_{i}^{0} \), where \( b_{i}^{0} = (1/m) \sum_{j=0}^{m-1} x_{i+j} \). Compute the similarity degree \( K_{ij}^m = \mu(d_{ij}^m, r) \), where \( d_{ij}^m \) is the maximum absolute difference of the components of \( b_i \) and \( b_j \), and \( \mu \) is a fuzzy function. Let 
$B_i^{m}(r) = ({1}/({N-m-1})) \sum_{\substack{j=1 \\ j\neq i}}^{N-m} K_{ij}^m$
and 
$B^{m}(r) = ({1}/({N-m})) \sum_{i=1}^{N-m} B_i^{m}(r),$
then 
$H_{FE}(m,r) = \lim_{N\to \infty} \left[\ln(B^{m}(r)) - \ln(B^{m+1}(r))\right]$.
 &  $H_{FE}(m,r)$ determines the similarity between two vectors based on their shapes rather than their absolute coordinates. Instead of the traditional Heaviside function, the fuzzy membership function is used. It shows greater consistency and less dependence on data length.
 \\

 $2$ & \citep{valencia2009refinedmultiscaleentropy} & Refined Multiscale Entropy 
 & Replace the finite impulse response filter with a low-pass butterworth filter with a squared magnitude of its frequency response is $|H(e^{j2\pi g})|^2=1/(1+{(g/g_c)}^{2n})$, where $n$ is filter order, and $g_c$ is the cutoff frequency. 
 & Continuously update the tolerable distance constant $r$ as a proportion of the standard deviation of the filtered series after coarse-graining, and then compute the sample entropy.  
 & To minimize the dependence of the computed entropy on the reduced variance.  
 \\
 $3$ & \citep{wu2013timecompositemultiscaleentropy} & Composite Multiscale Entropy 
 & The coarse-grained of the time sequence $\{x_i\}$ for a scale factor $d$ is given by $\{y_{l,j}^d\}_j$, where $y_{l,j}^d=\frac{1}{d}\sum_{i=(j-1)d+l}^{jd+l-1}x_i$, $1\leq j \leq N/d$ and $1\leq l \leq d$.
 & For each scale parameter \( d \), compute the corresponding sample entropy $H_{SE}(m,r,N/d)$ of the coarse-grained sequence $\{y_{l,j}^d\}_j$. The composite multiscale entropy is the average value $\frac{1}{d}\sum_{l=1}^{d}H_{SE}(m,r,N/d)$.
 & It provides improved results for short time series with high variance, particularly when the scale parameter in coarse-graining increases, leading to a rise in the variance of the corresponding time sequence.
 \\
 $4$ & \citep{wu2013modifiedmultiscaleentropy} & Modified Multiscale Entropy 
 & The coarse-graining procedure is replaced with a moving-average method, defined by $y_j^{d}=\frac{1}{d}\sum_{i=j}^{j+d-1}x_i$, where $1\leq j \leq N-d+1$.
 & Calculate the sample entropy of the time sequence $\{y_j^{d}\}$ using Equation \ref{sampleentropywiththreeparameter}.
 & It addresses the issues of inaccurate estimation and undefined entropy values in short-time sequences.
 \\

\bottomrule
\end{tabular}
\caption{Multiscale Entropy-Based generalized entropies: definitions and key properties}
\label{Table of generalized entropy based on multiscale entropy}
\end{sidewaystable}

The complexity measures discussed so far disregard the order of the time sequence and require a substantial amount of data to yield a meaningful entropy estimate \citep{riedl2013practicalpermutationentropyproperties}. Bandt and Pompe\citep{bandt2002permutationentropy} proposed \textbf{permutation entropy} as a simple and robust measure of complexity based on entropy functions and symbolic dynamics. Given a time sequence \(\{x_1,x_2,...,x_{N}\}\) and an embedding dimension \(E (\geq 2)\), we can construct \((N-E)\) vectors of consecutive values, each of length \(E\), given by
\begin{equation}
    v_i=(x_{i+0},x_{i+1},...,x_{i+E-1}) \ \ \ i=1,2,...,N-E+1.
\end{equation}
It is well-known that a vector of length \( E \), consisting of the elements \( \{0, 1, 2, \ldots, E-1\} \), can be arranged in \( E! \) different ways. Calculate the permutation \( \Tilde{q} = (q_0, q_1, \ldots, q_{E-1}) \) of \( (0, 1, 2, \ldots, E-1) \) such that 
\begin{equation}
    x_{i+q_0}\leq  x_{i+q_1} \leq  x_{i+q_2} \cdots \leq  x_{i+q_{E-1}},
\end{equation}
for each $v_i$, where $i=1,2,...,N-E+1$. Now, the probability of permutation $\Tilde{q}$ is estimated as 
\begin{equation}
    p(\Tilde{q})=\frac{\# \{i; i\leq N-E+1 \textit{ and $(x_{i+0},x_{i+1},...,x_{i+E-1})$ has the permutation $\Tilde{q}$} \}}{N-E+1}.
\end{equation}
Thus, the permutation entropy \citep{bandt2002permutationentropy} of order $E$ is given by 
\begin{equation}\label{permutationentropy}
    H_{PE}(E)=-\sum p(\Tilde{q}) \ln (p(\Tilde{q})),
\end{equation}
where summation is over all $E!$ permutations. When values are equal, they are ordered based on their occurrence time. It satisfies the inequality \(0 \leq H_{PE}(E) \leq \ln(E!)\). If the time sequence increases or decreases, then \(H_{PE}(E) = 0\) is achieved. Conversely, we get \(H_{PE}(E) = \ln(E!)\) for a completely random sequence (i.i.d. uniform random variables). The more frequently used measure is the normalized permutation entropy, defined by 
\begin{equation}\label{normalizedpermutationentropy}
    H_{PE}^{N}(E)=-\frac{1}{\ln(E!)}\sum p(\Tilde{q}) \ln (p(\Tilde{q})).
\end{equation}
It ranges between $0$ and $1$, where $0$ signifies completely predictable dynamics and $1$ represents completely stochastic dynamics. The commonly chosen embedding dimension $E$ ranges from $3$ to $7$. If the time sequence is independently and identically distributed, the statistic $2([N-E+1][\ln(E!)][1-H_{PE}^{N}(E)])$ asymptotically follows $\chi_{E!-1}^{2}$ \citep{huang2022permutationentropyproperties2}.

Permutation entropy relies entirely on the probability estimated from permutation patterns, which results in the loss of information related to amplitude. Consequently, permutation patterns with higher and lower amplitude values should not contribute equally to the permutation entropy. \textbf{Weighted permutation entropy} offers an effective method for assigning appropriate weights to pattern changes \citep{fadlallah2013weightedpermutationentropy}. In this approach, the frequency is estimated for each temporal pattern \(\Tilde{q}_k\) as
\begin{equation}
    p_w(\Tilde{q}_k)=\frac{w_k\#\{k; k\leq N-E+1 \textit{ and $(x_{k+0},x_{k+1},...,x_{k+E-1})$ has the permutation $\Tilde{q}_k$} \}}{\sum_{k}w_k\#\{k; k\leq N-E+1 \textit{ and $(x_{k+0},x_{k+1},...,x_{k+E-1})$ has the permutation $\Tilde{q}_k$} \}}.
\end{equation}
The weighted permutation entropy is defined as 
\begin{equation} \label{weightedpermutationentropy}
    H_{WPE}(E)=-\sum_{k} p_w(\Tilde{q}_k) \ln (p_w(\Tilde{q}_k)).
\end{equation}
It can be seen that when \( w_k = c \) for every \( k \), where \( c \) is a positive real constant, $H_{WPE}(E)=H_{PE}(E)$. The function in (\ref{weightedpermutationentropy}) retains most of the properties of permutation entropy and, under affine linear transformations remains invariant. The assignment of weights depends on the specific datasets used and the domain knowledge. $H_{WPE}(E)$ \citep{yin2016weightedpermutationentropyproperties} is adept at identifying sharp shifts in the signal and effectively distinguishes amplitude variations between identical ordinal patterns by assigning higher(lower) complexity to segments that are influenced by noise(exhibit regularity).

Analyzing the entire time sequence with precise amplitudes is computationally intensive and makes pattern interpretation challenging. In the \textbf{Dispersion entropy} algorithm \citep{rostaghi2016dispersionentropy}, the entropy measure is calculated by dividing the time sequence into a finite number of accessible classes. For a time sequence \( X=\{x_1,x_2,...,x_N\} \), the computation of dispersion entropy involves the following steps:
\begin{enumerate}
    \item Let \( m \) be a positive integer. Normalize each point in \( X \) to ensure \( 1 \leq x_i \leq m \) for all \( i \), using any linear or non-linear method.
    \item Divide the point sequence into \( m \) classes by multiplying each \( x_i \) by \( m \), adding 0.5, and then taking the nearest integer. Thus, 
    \begin{equation}
        y_i^m= \text{Nearest integer of } \{ mx_i+0.5 \}.
    \end{equation}
    \item Let the embedding dimension be \( k \) and the time delay \( t \). Compute \( N - (k-1)t \) time sequences 
    \begin{equation}
        w_{i}^{k,m}=\{ y_{i}^m,y_{i+t}^m,...,y_{i+(k-1)t}^m \}.
    \end{equation}
    \item Each \( y_{i}^m \) is an integer ranging from 0 to \( m \), resulting in \( m^k \) possible dispersion patterns for \( w_{j}^{k,m} \). Assign dispersion pattern ${\eta}_{u_0u_1...u_{k-1}}$ to $w_{j}^{k,m}$ as $u_0=y_{i}^m$, $u_1=y_{i+t}^m$, ..., $u_{k-1}=y_{i+(k-1)t}^m$.
    \item The relative frequency of pattern ${\eta}_{u_0u_1...u_{k-1}}$ is given by 
    \begin{equation}
        p({\eta}_{u_0u_1...u_{k-1}})=\frac{\#\{ j| j\leq N-(k-1)t, w_{j}^{k,m} \text{ \ \ has the pattern \ \ } {\eta}_{u_0u_1...u_{k-1}}  \}}{ N-(k-1)t }.
    \end{equation}
    \item Given the embedding dimension \( k \), time delay \( t \), and a number of classes \( m \), the dispersion entropy is defined as 
    \begin{equation}
        H_{DE}(X,k,t,m)=-\sum_{{\eta}_{u_0u_1...u_{k-1}}}p({\eta}_{u_0u_1...u_{k-1}})\ln (p({\eta}_{u_0u_1...u_{k-1}})),
    \end{equation}
    and the normalized dispersion entropy is given by 
    \begin{equation}
        H_{DE}^{N}(X,k,t,m)=-\frac{1}{\ln(m^k)}\sum_{{\eta}_{u_0u_1...u_{k-1}}}p({\eta}_{u_0u_1...u_{k-1}})\ln (p({\eta}_{u_0u_1...u_{k-1}})).
    \end{equation}
\end{enumerate}
There are several other measures to capture the uncertainty inherent in time-domain data, such as multiscale permutation entropy \citep{aziz2005multiscalepermutationentropy} and multivariate multiscale entropy \citep{ahmed2011multivariatemultiscalentropy}. We now describe entropy measures designed to quantify vagueness in the data, which are based on fuzzy theory \citep{zadeh1965fuzzysettheory}. This theory is adequate for modelling specific types of uncertainty, facilitating the determination of approximate solutions, and easing duality constraints.

\subsection{Fuzzy entropy measures}
According to the classical set theory, an element either belongs to a set or not, thus allowing for only a binary classification of membership. Fuzzy set theory \citep{zimmermann2011fuzzysettheorybook} utilises a membership function to quantify the degree of inclusion of an element within a set. Specifically, for any given set \( X \), a membership function is defined as a mapping from \( X \) to a subset of non-negative real numbers, possessing a finite supremum, thereby providing a measure of set membership for every element of $X$. The collection of ordered pairs  $(x, \mu_{\Tilde{X}}(x))$ for all $x\in X$, where $\mu_{\Tilde{X}}$ is the membership function, is defined as a fuzzy set $\Tilde{X}$.

Let \( Y \) be a random variable with a sample space \(S=\{y_1, y_2, \dots, y_n\}\) and a corresponding probability distribution \(\Tilde{P}=\{p_1, p_2, \dots, p_n\}\). Zadeh's \textbf{Fuzzy entropy} \citep{zadeh1968probabilityfuzzyentropyzadeh} associated with a fuzzy set \( \Tilde{S} \) is defined as 
\begin{equation} \label{fuzzyentropy}
    H_{FE}(\Tilde{S})=-\sum_{i=1}^{n}\mu_{\Tilde{S}}(y_i)p_i \ln(p_i),
\end{equation}
where $\Tilde{S}$ is a fuzzy set with membership function $\mu_{\Tilde{S}}$ on $S$.  $H_{FE}(\Tilde{S})$ quantifies the uncertainty inherent in the elements of the fuzzy set \( \Tilde{S} \). It is always non-negative. Note that if \( \Tilde{S} \) is non-fuzzy, \( H_{FE} \) does not reduce to the entropy of the distribution \( \Tilde{P} \), except when \( \Tilde{S} \) is the entire sample space. Consider the independent random variables \( Y_1 \) and \( Y_2 \) with probability distributions \( P=\{p_1, p_2,..., p_n\} \) and \( Q=\{q_1, q_2,...,q_m\} \), respectively. The fuzzy entropy of their joint probability distribution \( PQ={\{p_iq_j\}}_{i=1,j=1}^{n,m} \) is written as 
\begin{equation}
    H_{FE}(\Tilde{S}_1\Tilde{S}_2)=P^{\Tilde{S}_1} H_{FE}(\Tilde{S}_1)+ P^{\Tilde{S}_2} H_{FE}(\Tilde{S}_2), 
\end{equation}
where
\begin{equation}
    P^{\Tilde{S}_1}=\sum_{i=1}^{n}\mu_{\Tilde{S}_1}(y^1_i)p_i,
\end{equation}
\begin{equation}
    P^{\Tilde{S}_2}=\sum_{i=1}^{m}\mu_{\Tilde{S}_2}(y^2_i)q_i
\end{equation}
and $\Tilde{S}_1$ and $\Tilde{S}_2$ are the fuzzy sets associated to distributions $P$ and $Q$, respectively. 

Membership function is one of the key factor influencing the entropy of a fuzzy set. To address this effectively, three conditions are considered in \citep{de1993definitionmembershipfuzzyentropy}, which are used to construct complexity measures. Let \( \mathcal{L} \) be a lattice comprising maps from \( X \) to \([0,1]\), and \( \Psi \) a functional defined on \( \mathcal{L} \). The conditions proposed by Deluca \cite{de1993definitionmembershipfuzzyentropy} considered in the construction of the entropy function are- \( \Psi \) equals $0$ when \( f\in \mathcal{L} \) is either $0$ or $1$, reaches its maximum value when \( f = {1}/{2} \), and satisfies \( \Psi(f) \geq \Psi(g) \) if \( g \) is a sharpened version of \( f \), that is, for \( f(x) \geq {1}/{2} \) and \( f(x) \leq {1}/{2} \), \( g(x) \geq f(x) \) and \( g(x) \leq f(x) \), respectively. If \( \mathcal{L} \) is a finite lattice and $f\in \mathcal{L}$, then \textbf{Deluca and Termini fuzzy entropy} is defined as 
\begin{equation} \label{DTfuzzyentropy}
    H_{DTFE}(f)=-\mathcal{K}\sum_{i=1}^{n}\left( f(x_i) \ln(f(x_i))+(1-f(x_i))\ln(1-f(x_i)) \right),
\end{equation}
here $n$ is the size of $X(=\{x_1,...,x_n\})$ and $\mathcal{K}$ a positive constant. The functional in (\ref{DTfuzzyentropy}) satisfies all three conditions outlined above. $H_{DTFE}$ is a non-negative function on $\mathcal{L}$. We can write the functional \( H_{DTFE}(f) \) as \begin{equation}
    H_{DTFE}(f)=-\mathcal{K}\left(\sum_{i=1}^{n}h(f(x_i))+\sum_{i=1}^{n}h(1-f(x_i))\right),
\end{equation}
where $h(f(x_i))=f(x_i) \ln(f(x_i))$. Let \( f \) and \( g \) be two membership functions for a set \( X \), with their direct product defined as \( f*g(x,y) = f(x).g(y) \). If the power of membership function is given by \( F = \sum_{i} f(x_i) \) and \( G = \sum_{i} g(y_i) \), then we have the relation 
\begin{equation}
    h(f*g)=G.h(f)+F.h(g).
\end{equation}

Bruce \citep{ebanks1983measuresbrucefuzzyentropy} provides an axiomatic definition of the entropy function. Given a set \( X \) and a lattice \( \mathcal{L}=[0,1]^X \), \textbf{Bruce's fuzzy entropy} is defined as 
\begin{equation}\label{brucefuzzyentropy}
    H_{BFE}(f_X)=\sum_{x\in X}f_X(x)(1-f_X(x)), \text{ \ \ $f_X\in \mathcal{L}$.}
\end{equation}
The function in equation (\ref{brucefuzzyentropy}) is unique if and only if the following conditions hold.
\begin{enumerate}
    \item Sharpness: if $f_X\in\{0,1\}$ then $H_{BFE}(f_X)=0$.
    \item Maximality: $H_{BFE}(f_X)$ is maximum for $f_X\equiv 1/2$.
    \item Resolution: if $f_X^*$ is the sharpened version of $f_X$ then $H_{BFE}(f_X)\geq H_{BFE}(f_X^*)$.
    \item Symmetry: $H_{BFE}(f_X)=H_{BFE}(1-f_X)$.
    \item Valuation: For every $f_X,g_X \in \mathcal{L}$,
    \begin{equation}
        H_{BFE}(\max\{f_X,g_X\})+H_{BFE}(\min\{f_X,g_X\})=H_{BFE}(f_X)+H_{BFE}(g_X).
    \end{equation}
    \item Generalized additivity: Let $F_X$ and $G_Y$ be the powers of membership functions $f_X$ and $g_Y$, respectively. There exist functions $\phi, \psi:[0,\infty)\to [0,\infty)$ such that for finite sets $X$ and $Y$, we have
    \begin{equation}
        H_{BFE}(f_X*g_Y)=\phi(G_Y)H_{BFE}(f_X)+\psi(F_X)H_{BFE}(g_Y).
    \end{equation}
\end{enumerate}

Let \( X = \{x_1, x_2, \dots, x_n\} \) be a fuzzy set. Each subset of \( X \) can be represented as a bit vector by assigning $1$ if \( x_i \) is in the subset and $0$ if \( x_i \) is not; for example, if \( X = \{x_1, x_2, x_3\} \) and \( A = \{x_2\} \), then \( A = (0, 1, 0) \). Let $A=\{y_i\}$ be a fuzzy message or subset of $X$. The farthest non-fuzzy message $A^{F}$ is defined by assigning $1$ if $m_X(y_i) \leq 0.5$ and $0$ if $m_X(y_i) \geq 0.5$; conversely, the nearest non-fuzzy message $A^{N}$ is defined by assigning $0$ if $m_X(y_i) \leq 0.5$ and $1$ if $m_X(y_i) \geq 0.5$ \citep{fuzzyentropysurvey}. For example, if $A = \{0.1, 0.8, 0.2, 0.9, 0.5\}$, then $A^{F} = \{1, 0, 1, 0, 0\}$ (or equivalently $(1, 0, 1, 0, 1)$), and $A^{N} = \{0, 1, 0, 1, 0\}$ (or $(0, 1, 0, 1, 1)$). The \( l^p \) norm between the fuzzy messages \( A \) and \( B \) is given by 
\begin{equation}\label{fuzzydistance}
    l^p(A,B)={\left( \sum_{i}{|m_{A}(x_i)-m_{B}(x_i)|}^{p} \right)}^{1/p}, \textit{\ \ where \ \ } p\geq 1.
\end{equation}
Then, the geometry-based \textbf{Kosko fuzzy entropy} \citep{koskoentropyfuzzy} is a function \( H_{KFE}: P \to [0,1], \) defined by 
\begin{equation}\label{koskofuzzyentropy}
    H_{KFE}(X)=\frac{l^p(A,A^{N})}{l^p(A,A^{F})},
\end{equation}
where $p\geq1$ and $X$ is a fuzzy set. It satisfies the sharpness, maximality, resolution, and symmetry properties described above.
One major drawback of Kosko fuzzy entropy is its dependence on system resources. The lack of high-performance machines can result in significant computational overhead, as the fuzzy entropy calculation for a fuzzy set of $100$ elements requires constructing a $100$-dimensional hypercube.

The \textbf{Pal fuzzy entropy} \citep{palfuzzyentropy1989object} of a fuzzy set $A=\{x_i\}$ is defined by 
\begin{equation}
    H_{PFE}(A)=\frac{1}{n}\sum_{j=1}^{n}\left( m_A(x_j)e^{(1-m_A(x_j))}+(1-m_A(x_j))e^{m_A(x_j)} \right),
\end{equation}
where, $m_A$ is a membership function and $n$ is the size of the set $A$. It reflects the average level of uncertainty or ambiguity in determining whether an element belongs to set A or not. It satisfies several fundamental properties, including sharpness, maximality, resolution, and symmetry. 

Fuzzy set theory is a generalization of the classical set theory. Consequently, the entropy of fuzzy sets should also be a generalization. To address this, Pal et al. \citep{palHigherorderfuzzyentropy1992higher} introduced two definitions of fuzzy entropy to accommodate different scenarios. One is higher-order fuzzy entropy. In this context, suppose we want to measure uncertainty among the total \(n\) elements of a fuzzy set \(X = \{x_1, x_2, \dots, x_n\}\) possessing some property \(K\). The higher-order fuzzy entropy measures the average uncertainty corresponding to the possession of property \(K\) among a subset \(Y\) of \(X\) containing \(r\) elements. The \textbf{Higher-order Pal fuzzy entropy} of order $r$ for a fuzzy set $A$ is defined as 
\begin{equation}
    H_{HPFE}^r(A)=\left(1/\left( \Mycomb[n]{r} \right) \right)\sum_{i=1}^{\Mycomb[n]{r}}\left( m_A(S_i^r) \exp(1-m_A(S_i^r))+(1-m_A(S_i^r))\exp{m_A(S_i^r)} \right),
\end{equation}
where $m_A$ is a membership function, $S_i^r$ represents set of $r$ elements from $A$ and  
\begin{equation}
    m_A(S_i^r)=\min_{z\in {S_i^r} } \{ m_A(z) \}.
\end{equation}
It possesses the properties of sharpness, maximality, and resolution. The symmetry property for $H_{HPFE}^r(A)$ does not always hold in general. $H_{HPFE}^r(A)\geq H_{HPFE}^{r+1}(A)$, when $m_A$ is in $[0, 0.5]$, and $H_{HPFE}^r(A)\leq H_{HPFE}^{r+1}(A)$, otherwise. $H_{HPFE}^r(A)$ is a generalization of $H_{HPFE}(A)$ when $S_i^1$ is a singleton subset of $A$. The second type is the \textbf{Hybrid Pal fuzzy entropy} \citep{palHigherorderfuzzyentropy1992higher}. Let the binary symbols $0$ and $1$ occur with probability $p_0$ and $p_1(=1-p_0)$, respectively, and let the closeness of a symbol to $1$ be indicated by the membership function $m_A$. The hybrid Pal fuzzy entropy of fuzzy set (of symbols) $A$ is then defined by 
\begin{equation}
    H_{HbPFE}(A)=-p_0 \log(E_0)-p_1 \log(E_1),
\end{equation}
where,
\begin{equation}
E_0=\frac{1}{n}\sum_{i=1}^{n}(1-m_A(x_i))\exp(m_A(x_i)),
\end{equation}
represents the average likelihood of treating the received symbol as $0$ and
\begin{equation}
E_1=\frac{1}{n}\sum_{i=1}^{n}(m_A(x_i))\exp(1-m_A(x_i)),
\end{equation}
represents the average likelihood of treating the received symbol as $1$. $H_{HbPFE}$ reduces to discrete Shannon entropy in the absence of fuzziness. If \( m_A(x_i) = 0.5 \) for all $i$, then $H_{HbPFE}$ will be a constant function.

Let \( A \) be a fuzzy set, and $m_A$ be a corresponding membership function. Let \( A^* \) be the most fuzzy set (i.e. $m_A(x)=0.5$ for each $x$) corresponding to \( A \). Then, the \textbf{Bhandari and Pal fuzzy entropy} \citep{bhandaripalfuzzyentropy1993some} is defined by 
\begin{equation}
    H_{BPFE}(A)=\sum_{i=1}^{n}\left[ m_A(x_i)\ln\left(\frac{m_A(x_i)}{m_{A^*}(x_i)} \right)+(1-m_A(x_i))\ln\left(\frac{1-m_A(x_i)}{1-m_{A^*}(x_i)} \right) \right].
\end{equation}
It is related to the Deluca and Termini fuzzy entropy from equation (\ref{DTfuzzyentropy}), expressed as \( H_{DTFE}(A) = 1-\mathcal{K}H_{BPFE}(A) \). Given \( A^*=B \) as any other fuzzy set, the function corresponding to $H_{BPFE}(A)$ will be
\begin{equation}
    D_{BPFD}(A,B)=I_{BPFD}(A,B)+I_{BPFD}(B,A),
\end{equation}
where
\begin{equation}
I_{BPFD}(A,B)=\sum_{i=1}^{n}\left[ m_A(x_i)\ln\left(\frac{m_A(x_i)}{m_{B}(x_i)} \right)+(1-m_A(x_i))\ln\left(\frac{1-m_A(x_i)}{1-m_{B}(x_i)} \right) \right],
\end{equation}
which is commonly utilized as a measure of fuzzy divergence. It is always positive and equals zero if and only if \( A = B \). Additionally, \( D_{BPFD}(A,B) \) is a symmetric function. It possesses the following properties: 
\begin{enumerate}
    \item $D_{BPFD}(A\cup B,A\cap B)=D_{BPFD}(A,B)$.
    \item $D_{BPFD}(A\cup B,C)\leq D_{BPFD}(A,C) + D_{BPFD}(B,C)$, where $C$ is a fuzzy set.
    \item $D_{BPFD}(A,B)$ is maximum iff the farthest nonfuzzy subset of $A$ is $B$.
\end{enumerate}
Bhandari and Pal\citep{bhandaripalfuzzyentropy1993some} also defined a \textbf{fuzzy entropy of order \(\alpha\)}, inspired by Renyi's entropy, as
\begin{equation}\label{renyitypeBPFEfuzzyentropy}
    H_{aBPFE}^{\alpha}(A)=\frac{1}{c(1-\alpha)}\sum_{i=1}^{n}\ln\left({m_A(x_i)}^{\alpha}+{(1-m_A(x_i))}^{\alpha} \right), \ \ \ \alpha >0(\neq 1),
\end{equation}
where $c$ is the normalizing constant. Some properties are:
\begin{enumerate}
    \item $H_{aBPFE}^{\alpha}(A)=0$ iff $A$ is a nonfuzzy set.
    \item $H_{aBPFE}^{\alpha}(A)$ is maximum iff \(A\) is the most fuzzy set.
    \item $H_{aBPFE}^{\alpha}$ decreases for sharpened set.
    \item $H_{aBPFE}^{\alpha}(A\cup B)+ H_{aBPFE}^{\alpha}(A\cap B)=H_{aBPFE}^{\alpha}(A)+H_{aBPFE}^{\alpha}(B)$.
\end{enumerate}

To analyze and utilize fuzziness efficiently, Hooda \citep{hoodaparametricfuzzyentropy2004generalized} proposed one- and two-parametric \textbf{Hooda fuzzy entropy} measures, using the one- and two-parametric entropy measures provided in \citep{sharmaandmittalentropy1975new}, defined for a fuzzy set $A=\{x_1, x_2,..., x_n\}$ by 
\begin{equation}
    H_{HFE}^{\beta}(A)=\frac{1}{1-\beta}\left[ 2^{(\beta-1)\sum_{i=1}^{n}m_{A}(x_i)\log(m_{A}(x_i))+(1-m_{A}(x_i))\log(1-m_{A}(x_i))}-1 \right],
\end{equation}
where $m_{A}$ is a membership function, $\beta>0(\neq 1)$ and
\begin{equation}
    H_{HFE}^{\beta,\alpha}(A)=\frac{1}{1-\beta}\sum_{i=1}^{n}\left[{\left(m_{A}^{\alpha}(x_i)+{(1-m_{A}(x_i))}^{\alpha}\right)}^{\frac{\beta-1}{\alpha-1}}-1 \right],
\end{equation}
where $\alpha\neq \beta$, $\alpha, \beta >0$ and $\alpha\neq 1$. It is straightforward to observe that for \( \mathcal{K}=1 \), 

\begin{equation}
H_{HFE}^{\beta}(A) = h_{\beta}(H_{DTFE}(A)), 
\end{equation}
where $H_{DTFE}(A)$ is the Deluca and Termini fuzzy entropy as in equation (\ref{DTfuzzyentropy}) and the function 
\begin{equation}
h_{\beta}(x) = {(1-\beta)}^{-1}\left(2^{(1-\beta)x}-1\right),
\end{equation}
for $x\geq0$. $H_{HFE}^{\beta}$ satisfies the properties of sharpness, maximality, resolution, and symmetry. Note that for \(\beta = 1\), with appropriate adjustments, $H_{HFE}^{1}$ will be equal to \(H_{DTFE}\) as in equation (\ref{DTfuzzyentropy}), for $\mathcal{K}=1$. Additionally, for $c=1$, the relation 
\begin{equation}
H_{HFE}^{\beta,\alpha}(A) = \sum_{i=1}^{n}h_{\beta}\left( W^{\alpha}(m_A(x_i))  \right),
\end{equation}
where 
\begin{equation}
W^{\alpha}(m_A(x_i)) =\frac{1}{1-\alpha}\left( \log(m_A^{\alpha}(x_i)+{(1-m_A(x_i))}^{\alpha}) \right),
\end{equation}
can be proved. $H_{HFE}^{\beta,\alpha}(A)$ is a Convex downward function, attaining its maximum if and only if A is maximally fuzzy, that is, when $m_{A}(x_i)=0.5$ for all $i$. It also satisfies the maximality, resolution and symmetry properties. For $\beta=1$, the relation 
\begin{equation}
    H_{HFE}^{1,\alpha}(A)=H_{aBPFE}^{\alpha}(A),
\end{equation}
with $c=1$ from equation (\ref{renyitypeBPFEfuzzyentropy}), holds.

One generalization of a fuzzy set $X=\{x\}$ is the intuitionistic fuzzy set \citep{atanassov1999intuitionisticfuzzyset}. Let $m_{X}:X\to [0,1]$ and $n_{X}:X\to [0,1]$ be the function of membership and non-membership degrees on $X$, respectively. Similar to the definition of a fuzzy set as $\{x, m_X(x)\}$ for $x \in X$, the collection $\{x, m_{X}(x), n_{X}(x)\}$ for $x \in X$ constitutes an intuitionistic fuzzy set if $0\leq m_{X}(x)+n_{X}(x)\leq 1$ for all $x\in X$. We first define the hesitation margin ($q_X$) before defining entropy for intuitionistic fuzzy sets. For a given \( x \) in \( X \), it is defined by 
\begin{equation}
    q_X(x)=1-m_{X}(x)-n_{X}(x).
\end{equation}
Also, the normalized Hamming distance \citep{szmidt2000distancebetweenelementsoffuzzyset} between two intuitionistic fuzzy sets, $A$ and $B$, is given by 
\begin{equation}
    d(A,B)=\frac{1}{2n}\sum_{i=1}^{n}\left( |m_{A}(x_i)-m_{B}(x_i)|+|n_{A}(x_i)-n_{B}(x_i)|+|q_{A}(x_i)-q_{B}(x_i)| \right).
\end{equation}
Let \( S \) denotes the elements that fully belong to the set, that is \( m_X = 1\) and \( n_X=0 \), and \( T \) denotes the elements that do not fully belong to the set, that is \( n_X = 1\) and \( m_X=0 \). Then the intuitionistic fuzzy entropy of an element, \( x \in X \), is given by 
\begin{equation}
    ent(x)=\frac{d_n}{d_f},
\end{equation}
where $d_n$ and $d_f$ are the smallest and largest hamming distances of $x$ from the elements of the set $S$ and $T$ respectively. It quantifies the information needed to determine whether an element \( x \), characterized by \( (m_X, n_X, q_X) \), either fully belongs to or does not belong to the set. Thus, for a $n$ element set $X$, the \textbf{intuitionistic fuzzy entropy} \citep{szmidt2014measureintuitionisticfuzzyentropy} is defined by 
\begin{equation}
    Ent(X)=\frac{1}{n}\sum_{i=1}^{n}ent(x_i).
\end{equation}
For all values of the hesitation margin $q_{X}$, \( Ent(X) \) reaches its maximum when \( m_X = n_X \). It is always a positive measure that increases from point \( S \) to the centre between \( S \) and \( T \) and then decreases until it reaches to \( T \). 

\subsection{Fractional order entropy measures}
Now, we summarize the fractional-order entropy measures. Fractional calculus provides a generalized framework for modelling complex systems with memory and hereditary properties, enabling more accurate descriptions of processes such as anomalous diffusion and viscoelastic behaviour. For a concise introduction, refer to \citep{oliveira2014reviewfractionalcalculus}.

It is known that Shannon and Tsallis entropies can be expressed, respectively, as a limit of the standard and Jackson q-derivatives of a function of probabilities. Based on this idea, a one-parameter \textbf{Akimoto-Suzuki fractional-order entropy} \citep{ASfractionalorderentropy} is defined by 
\begin{equation}
    H_{ASFOE}^{\delta}=-\lim_{x \to 1}\sum_{i}\dv[\delta]{x} e^{x \ln p_i},
\end{equation}
where $\dv[\delta]{f(x)}{x}=\leftidx{^a}D_{RL}^{\delta}f(x)$ is the Riemann-Liouville derivative with $a=0$, which has the left(right)-hand side derivative given by

\begin{equation}
    \leftidx{^a}D_{RL}^{\delta}f(x)=\frac{1}{\Gamma(n-\delta)}\dv[n]{x}\int_{a(x)}^{x(a)}\frac{f(z)}{{(x-z)}^{\delta-n+1}}dz, \ \ \ x\geq(\leq) a,
\end{equation}
here $n$ is such that $n-1<\delta<n$, $n\in \mathbb{N}$. $H_{ASFOE}^{\delta}$ possesses the properties of concavity, non-extensivity, and positivity. As the parameter \( \delta \) approaches $1$, \( H_{ASFOE}^{\delta} \) converges to the Shannon entropy as in equation (\ref{shannonentropy}). It is an increasing function with respect to the sample size. 

Next, one-parameter \textbf{Ubriaco fractional order entropy} \citep{ubriaco2009entropiesUfractionalorderentropy} is defined as 
\begin{equation}
    H_{UFOE}^{\delta}=\lim_{x\to -1}\dv{x}\left( \leftidx{^{-\infty}}D_{RL}^{1-\delta}\sum_{i}e^{-x\ln(p_i)} \right),
\end{equation}
where $\leftidx{^{-\infty}}D_{RL}^{1-\delta}$ is the left-hand side Riemann-Liouville derivative with $a\to -\infty$. To simplify this, we can express it as
\begin{equation}
    H_{UFOE}^{\delta}=\sum_{i}p_i{(-\ln p_i)}^{\delta}, \ \ \ 0\leq \delta \leq 1.
\end{equation}
$H_{UFOE}^{\delta}$ is concave, positive-definite, non-additive and satisfies the Lesche stability criteria\citep{leschestabilities} and thermodynamic stability\citep{thermodynamicstability} properties. As \( \delta \) approaches $1$, in the limiting case, $H_{UFOE}^{\delta}$ equals the Shannon entropy.

\noindent The information measure corresponding to Shannon entropy, 
\begin{equation}
\mathcal{I}(p_i) = -\ln(p_i), \ \ \ \forall i=1,2,...,n.
\end{equation}
Capitalizing on this relationship, Machado \citep{machadofractionalorderentropy} defined the generalized information measure by 
\begin{equation}\label{machadoinformationmeasure}
    \mathcal{I}^{\delta}(p_i)=\leftidx{^{a+}}D_{RL}^{\delta}\mathcal{I}(p_i)=-\frac{{p_i}^{-\delta}}{\Gamma(\delta+1)}\left( \ln (p_i)+\Psi(1)-\Psi(1-\delta) \right),
\end{equation}
where $\Psi$ is used for the digamma function and $\delta$ is a real number.
Further, the \textbf{Machado fractional order entropy} is defined by
\begin{equation}
    H_{MFOE}^{\delta}=\sum_{i}{p_i}\mathcal{I}^{\delta}(p_i).
\end{equation}
$H_{MFOE}^{\delta}$ can assume both positive and negative values. The one-parameter fractional entropy does not fully adhere to all Shannon entropy axioms, except when \( \delta = 0 \), where it converges to Shannon entropy. This aligns with the general principle that certain properties may be lost in the generalisation process.

The Karci fractional derivative \citep{lopes2020reviewoffractionalorderentropy} of a function \( h \) is defined by 
\begin{equation}
    D_{\delta}^{K}h(x)=\frac{\dv{x}[h(x)].{[h(x)]}^{\delta-1}}{x^{\delta-1}},
\end{equation}
where $[.]$ gives the integer part of the function and for some non-negative integer $n$, $n-1<\delta<n$. By applying the Karci fractional derivative to the differential form $H=\lim_{x\to -1}\sum_{i}p_i^{-x}$, of Shannon entropy, we obtain the \textbf{Karci fractional order entropy} \citep{karci2016karcifractionalorderentropy}, defined as 
\begin{equation}
    H_{KFOE}^{\delta}=\sum_{i}p_i.|{(-p_i)}^{\delta}\ln p_i|.
\end{equation}
It is always a positive real number. In the fractional order generalization of entropy based on Rényi entropy, Machado and Lopes\citep{machadorenyifractionalorderentropy} propose two definitions of entropy. The first type of \textbf{Machado and Lopes fractional order entropy} utilizes the Machado information measure, as given in equation (\ref{machadoinformationmeasure}), in the expected information form of Rényi's entropy 
\begin{equation}
    H^{\alpha}=\frac{1}{1-\alpha}\sum_{i}p_i. e^{(1-\alpha)\mathcal{I}(p_i)},
\end{equation}
with the final form represented as 
\begin{equation}
    H_{MLFOE1}^{\alpha,\delta}=\frac{1}{1-\alpha}\ln \left\{ \sum_{i} p_i \exp\left[ (\alpha-1).\frac{p_i^{-\delta}}{\Gamma(\delta+1)}\left(\ln(p_i)+\Psi(1)-\Psi(1-\delta)\right) \right]  \right\}.
\end{equation}
The second definition is obtained by modifying the Rényi entropy function (\ref{renyientropydiscrete}) and using (\ref{machadoinformationmeasure}), resulting in the measure given by 
\begin{equation}
    H_{MLFOE2}^{\alpha,\delta}=\frac{1}{n^{\frac{\delta}{\alpha}}}\frac{\alpha}{1-\alpha}\left[\frac{{\left(\frac{1}{n}\sum_{i}p_i^{\alpha}\right)}^{{-\delta}/{\alpha}}}{\Gamma(\delta+1)}\left( \frac{\ln(n)}{\alpha}+\ln\left({\left(\frac{1}{n}\sum_{i}p_i^{\alpha}\right)}^{{1}/{\alpha}}\right)+\Psi(1)-\Psi(1-\delta) \right)  \right],
\end{equation}
where $n$ is the sample size. The study highlights the flexibility and added degrees of freedom introduced by the two-parameter formulation in generalized fractional order entropy.

\subsection{Graph entropy measures}
 Let \( \mathcal{V} \) denote a set of vertices, and a set of edges \( \mathcal{E} \subseteq \{\{x, y\} \mid x, y \in V \text{ and } x \neq y\} \) consists of unordered pairs of distinct vertices then define an ordered pair \( \mathcal{G} = (\mathcal{V}, \mathcal{E}) \). In this context, \( \mathcal{G} \) is a graph, with each edge representing a connection between two vertices. A graph \( \mathcal{G} \) is connected if each pair of vertices in \( \mathcal{G} \) links by at least one path. When referring to a partition of a graph \(\mathcal{G}\), we specifically mean a partition of its vertex set. A bijective map \( \psi \) from a graph \( \mathcal{G} = (\mathcal{V}, \mathcal{E}) \) to itself is termed an automorphism on \( \mathcal{G} \) if it satisfies \( \psi(v_1 v_2) = \psi(v_1) \psi(v_2) \) where \( v_1 v_2, \, \psi(v_1) \psi(v_2) \in \mathcal{E} \). Under the composition of maps, the set of all such automorphisms, denoted \(\operatorname{Aut}(\mathcal{G})\), forms a group. Let \(\mathcal{G}\) be a graph, and let \(H \leq \text{Aut}(\mathcal{G})\) be a subgroup of automorphisms group of \(\mathcal{G}\). Two vertices \(v_1\) and \(v_2\) are said to be similar under \(H\) if there exists an automorphism in \(H\) that maps \(v_1\) to \(v_2\). The equivalence classes defined by this relation are called the orbits of the graph under \(H\). Two vertices from the same orbit of a graph are considered topologically equivalent. A detailed explanation of graph theory can be found in \citep{west2001graphtheoryintroduction}.

 The \textbf{Rashevsky graph entropy} \citep{1955rashevskygraphentropy} of a connected graph $\mathcal{G}=(\mathcal{V},\mathcal{E})$ is given as 
\begin{equation} \label{graph1ras}
H_{RGE}(\mathcal{G})=-\sum_{j=1}^{n^v}\frac{N^v_i}{|\mathcal{V}|}\ln\left(\frac{N^v_i}{|\mathcal{V}|}\right),
\end{equation}
where $N^v_i$ is the number of topological equivalent vertices in $i$th orbit, and $n^v$ is the number of orbits. It quantifies the structural complexity of a graph. It is based on the definition of Shannon entropy and, therefore, adheres to its properties such as positivity, expansibility and maximum for singleton orbits. 

Similarly, \textbf{Trucco graph entropy} \citep{Truccographentropy1956ANO} is defined, based on the edge automorphisms of a graph $\mathcal{G}$, as 
\begin{equation} \label{graph2ras}
    H_{TGE}(\mathcal{G})=-\sum_{i=1}^{n^E}\frac{N^E_i}{|\mathcal{E}|}\ln\left(\frac{N^E_i}{|\mathcal{E}|}\right),
\end{equation}
here $n^E$ is the number of orbits, and $N^E_i$ is the number of topological equivalent edges in each orbit. One distinction between (\ref{graph1ras}) and (\ref{graph2ras}) lies in their definitions: (\ref{graph1ras}) is based on a partition of the vertex set \( V \), whereas (\ref{graph2ras}) is defined on a partition of the edge set \( E \) through equivalent relation.

An entropy measure based on graph invariants, such as the number of edges and vertices, can yield the same value for structurally non-equivalent graphs. For instance, two non-isomorphic graphs may have identical \( H_{RGE} \) and \( H_{TGE} \) values. To overcome this, let  D\( = \left(d_{ij}\right) \) represents the distance matrix of a graph \( \mathcal{G} \), where \( d_{ij} \) denotes the distance units between vertices values \(1\leq i,j \leq |\mathcal{V}|\). For example, refer to graphs I and II and the corresponding matrices D(I) and D(II).

\begin{figure}[htbp]
    \begin{minipage}{0.5\textwidth}
        \centering
        \begin{tikzpicture}
            \foreach \i in {1,2,3,4} {
                \draw[thick] (\i-1+0.15,0) -- (\i-0.15,0); 
            }

            \foreach \i in {1,...,5} {
                \draw[thick] (\i-1,0) circle(0.15); 
                \node[above] at (\i-1,0.25) {\i}; 
            }
        \end{tikzpicture}
        \captionsetup{labelformat=empty} 
        \caption{I}
    \end{minipage}%
    \hfill
    \begin{minipage}{0.5\textwidth}
        \centering
        \begin{tikzpicture}
            \foreach \i in {1,2,3} {
                \draw[thick] (\i-1+0.15,0) -- (\i-0.15,0); 
            }

            \foreach \i in {1,...,4} {
                \draw[thick] (\i-1,0) circle(0.15); 
                \node[above] at (\i-1,0.25) {\i}; 
            }

            \draw[thick] (2,-0.15) -- (2,-0.85); 
            \draw[thick] (2,-1) circle(0.15); 
            \node[below] at (2,-1.25) {5}; 
        \end{tikzpicture}
        \captionsetup{labelformat=empty} 
        \caption{II}
    \end{minipage}
    \label{graphsexample}
\end{figure}

\begin{figure}[htbp]
    \begin{minipage}{0.5\textwidth}
        \centering
        \[
        \large
        \text{D(I)=} \quad
        \begin{pmatrix}
        0 & 1 & 2 & 3 & 4 \\
        1 & 0 & 1 & 2 & 3 \\
        2 & 1 & 0 & 1 & 2 \\
        3 & 2 & 1 & 0 & 1 \\
        4 & 3 & 2 & 1 & 0
        \end{pmatrix}
        \]
    \end{minipage}%
    \hfill
    \begin{minipage}{0.5\textwidth}
        \centering
        \[
        \large
        \text{D(II)=} \quad
        \begin{pmatrix}
        0 & 1 & 2 & 3 & 3 \\
        1 & 0 & 1 & 2 & 2 \\
        2 & 1 & 0 & 1 & 1 \\
        3 & 2 & 1 & 0 & 2 \\
        3 & 2 & 1 & 2 & 0
        \end{pmatrix}
        \]
    \end{minipage}
    \label{matrixforgraphsexample}
\end{figure}

Let \( n_k \) denote the number of \( k \) appears in the distance matrix where $1\leq k \leq |\mathcal{V}|-1$. The first type of \textbf{Bonchev and Trinajstić graph entropy} \citep{bonchev1977informationBTgraphentropy} is defined as 
\begin{equation}
    H_{BTGEa}(\mathcal{G})= {|\mathcal{V}|}^2 \ln({|\mathcal{V}|}^2)-|\mathcal{V}|\ln\left(|\mathcal{V}|\right)-\sum_{k=1}^{|\mathcal{V}|}n_k \ln\left(n_k\right).
\end{equation}
Since D is a symmetric matrix, the upper triangular portion of D is sufficient for computing \(  H_{BTGEa}(\mathcal{G}) \). Let $W^{\mathcal{G}}$ denote the Wiener number, capture information about the distribution of distances in a graph $\mathcal{G}$,
\begin{equation}
W^{\mathcal{G}}=\sum_{k=1}^{|\mathcal{V}|}k\frac{n_k}{2}. 
\end{equation}
Here \( H_{BTGEa} \) represents the regularity in a graph\citep{bonchev1977informationBTgraphentropy}; as the branching of the graph increases, \( H_{BTGEa} \) decreases. The Wiener number, \( W^{\mathcal{G}} \), quantifies the centrality of the graph, with a higher value of \( W^{\mathcal{G}} \), indicating a more complex and widely dispersed vertex structure. Based on this, the second type of \textbf{Bonchev and Trinajstić graph entropy} \citep{bonchev1977informationBTgraphentropy} is defined as 
\begin{equation}
   H_{BTGEb}(\mathcal{G})= W^{\mathcal{G}} \ln (W^{\mathcal{G}})-\sum_{k=1}^{|\mathcal{V}|} \frac{n_k}{2} k\ln (k).
\end{equation}
\( H_{BTGEb}(\mathcal{G}) \) is capable of distinguishing between graphs with varying Wiener numbers. It exhibits greater sensitivity compared to the Wiener number, as its value can vary not only with \( W^{\mathcal{G}} \) but also with changes in the distribution influenced by \( k \) and \( n_k \).

Measuring the uncertainty of subparts of a graph is also important; for instance, one may wish to calculate the uncertainty around a specific vertex. The number of distinct edges connected to \( v \) in a graph $\mathcal{G}$ is known as the degree $d(v)$ of a vertex \( v \).  Let
\begin{equation}
d = \sum_{i=1}^{m} n_i d(v_i),
\end{equation}
and
\begin{equation}
n =\sum_{i=1}^{m}n_i,
\end{equation}
where $n_i$ is the number of vertices of degree $d(v_i)$ and $m$ is the number of distinct possible degrees. Thus, we can obtain a probability distribution in the manner illustrated by the matrix \( A \) below.
\[
A=\begin{pmatrix}
d(v_1) & d(v_2) & \cdots & d(v_m) \\
n_1 & n_2 & \cdots & n_m \\
q_1 & q_2 & \cdots & q_m
\end{pmatrix}
\]
Here $q_k=\frac{d_k}{d}$ and $\sum_{i=1}^{m}n_iq_i=1$. The \textbf{Raychaudhury graph entropy} \citep{raychaudhury1984discriminationRRGRBgraphentropy} of graph \( \mathcal{G}=(\mathcal{V}, \mathcal{E}) \) calculates the degree of complexity and is given by 
\begin{equation}
\mathcal{I}^{d}_{RGEa}(\mathcal{G}) = -\sum_{k=1}^{m}n_k q_k\log_{2}(q_k).
\end{equation}
 The greatest distance of a vertex $v$ to any other vertex in the graph is known as the eccentricity $e(v)(=e)$ of a vertex $v$. The distance code for vertex \( v \) is expressed as \( v: 0^{1}, 1^{f_1}, 2^{f_2}, \dots, e^{f_e} \), where \( f_i \) represents the number of vertices at a distance of $i$ from \( v \). The distance frequency sequence is also given by \(\mathcal{D}(v)= (1, f_1, f_2, \dots, f_e) \). Let \( \mathcal{T}(v)=\sum_{j=1}^{e} j f_j \) represents the total distance from vertex \( v \). The probability distribution corresponding to the distance frequency sequence as the partition of $\mathcal{V}$ is given by matrix 
\[
B=\begin{pmatrix}
1 & f_1 & f_2 & \cdots & f_e \\
q_0^{'} & q_1^{'} & q_2^{'} & \cdots & q_e^{'}
\end{pmatrix}
\]
where \( q_0^{'} = {1}/{n} \) and \( q_i^{'} = {f_i}/{n} \) for $i=1,2,...,e$. The Raychaudhury vertex complexity \citep{raychaudhury1984discriminationRRGRBgraphentropy} of \( v \) is given as
\begin{equation}
h_{RGEb}(v) = \frac{1}{n}\log_2 (n)-\sum_{i=1}^{e}q_i^{'} \log_{2}(q_i^{'}),
\end{equation}
and the corresponding \textbf{Raychaudhury graph vertex complexity entropy} \citep{raychaudhury1984discriminationRRGRBgraphentropy} is defined as
\begin{equation}
H_{RGEb}(\mathcal{G}) = \frac{1}{n}\sum_{j=1}^{n}h_{RGEb}(v_j).
\end{equation}
Moreover, if the partition and probability scheme for \( \mathcal{T}(v) \) is given by matrix

\[
E=\begin{pmatrix}
1 & 2 & \cdots & e \\
f_1 & f_2 & \cdots & f_e \\
q_1^{''} & q_2^{''}  & \cdots & q_e^{''} 
\end{pmatrix},
\]
 where $q_i^{''}=i/\mathcal{T}(v)$ and $\sum_{i=1}^{e}f_i q_{i}^{''}=1$, then the corresponding vertex distance complexity measure for \( v \) is defined as 
\begin{equation}
h_{RGEc}(v) = -\sum_{i=1}^{e}f_i q_i^{''} \log_2 \left( q_i^{''} \right) , 
\end{equation}
and the \textbf{Raychaudhury graph distance complexity entropy} \citep{raychaudhury1984discriminationRRGRBgraphentropy} is given by
\begin{equation}
H_{RGEc}(\mathcal{G}) = \sum_{i=1}^{n}s_i h_{RGEc}(v_i),
\end{equation}
where $s_i={\mathcal{T}(v_i)}/{T}$ and $T=\sum_{i=1}^{n}{\mathcal{T}(v_i)}$. All three entropy functions satisfy properties similar to those of Shannon entropy, as their definitions are aligned with it.

An undirected graph is one in which all edges are bidirectional. Let \( \mathcal{G}=\left( \mathcal{V}, \mathcal{E} \right) \) be a finite graph. We now introduce the generalized graph entropy functions. Let \( f: S \to \mathcal{G} \) be an abstract information function for the graph \( \mathcal{G} \) corresponding to a given set $S$. The \textbf{Dehmer generalized graph entropy} \citep{dehmer2008informationdehmergraphentropy} for $\mathcal{G}$ is defined as 
\begin{equation}
    H_{DGGEa}(\mathcal{G})=-\sum_{i=1}^{|\mathcal{V}|}\frac{f(v_i)}{\sum_{j=1}^{|\mathcal{V}|}f(v_j)}\log\left( \frac{f(v_i)}{\sum_{j=1}^{|\mathcal{V}|}f(v_j)} \right).
\end{equation}
It provides a range of entropy functions corresponding to the chosen function \( f \) based on metrical graph properties. 

Dehmer also provides an algorithm for decomposing a graph into local information graphs. A sequence of vertices and edges of the graph connecting one vertex to another is a path in a graph. The path length is determined by the number of edges it includes. Consider an undirected connected graph \( \mathcal{G} \), for \( v_i \in \mathcal{V} \) and \( j = 1, 2, \ldots, \bar{m} \), let \( L^{\mathcal{G}}(j, v_i) \) denote the local information graph induced by the paths \( \mathcal{C}_1^j(v_i), \mathcal{C}_2^j(v_i), \ldots, \mathcal{C}_{k_j}^j(v_i) \). Let \( \mathcal{L}\left(L^{\mathcal{G}}(j, v_i) \right) = \sum_{w=1}^{k_j} l\left(\mathcal{C}_w^j(v_i)\right) \), where \( l\left(\mathcal{C}_k^j(v_i)\right) \) represents the length of \( \mathcal{C}_k^j(v_i) \) and $\bar{m}=\max_{v\in \mathcal{V}}e(v)$. Then, the information functional takes the form 
\begin{equation}
f(v_i) ={\gamma}^{{a_1}{\mathcal{L}\left(L^{\mathcal{G}}(1, v_i) \right)}+{a_2}{\mathcal{L}\left(L^{\mathcal{G}}(2, v_i) \right)}+\cdots+{a_{\bar{m}}}{\mathcal{L}\left(L^{\mathcal{G}}(\bar{m}, v_i) \right)}},
\end{equation}
where $\gamma>0$ and $a_i >0$ for all $i=1,2,\cdots,\bar{m}$.

Let \(\mathcal{G}=\left(\mathcal{V},\mathcal{E}\right)\) be a connected, undirected graph. To compute the generalized tree of height $\omega$, choose a vertex \(v_i\), use the algorithm to derive the generalized tree described in \citep{dehmer2008informationdehmergraphentropy} and apply this process to all the vertices in \(\mathcal{V}\). This yields the sequence of generalized trees \(GT(\mathcal{G}) = (H_1, H_2, \ldots, H_{|\mathcal{V}|})\). The vertex and edge \textbf{Dehmer graph entropy} based on the generalized tree are defined as 
\begin{equation}\label{dehmergraphentropyGT1}
H_{DGEb}^{v}(\mathcal{G}) = \sum_{i=1}^{|\mathcal{V}|} {h}_v(H_i)
\end{equation}
and
\begin{equation}\label{dehmergraphentropyGT2}
H_{DGEb}^{e}(\mathcal{G}) = \sum_{i=1}^{|\mathcal{V}|} {h}_e(H_i), 
\end{equation}
respectively. Here, 
\begin{equation}
    {h}_{v(e)}(H_i)=-\sum_{i=1}^{\omega}\frac{f^{v(e)}(l_i)}{\sum_{j=1}^{|\mathcal{V}|}f^{v(e)}(l_j)}\log\left( \frac{f^{v(e)}(l_i)}{\sum_{j=1}^{|\mathcal{V}|}f^{v(e)}(l_j)} \right),
\end{equation} 
and for some $\varrho>0,$
\begin{equation}
    f^{v(e)}(l_i)=\varrho^{|\mathcal{V}_i|(|\mathcal{E}_i|)}, \ \ \ \textit{ \ for each level $l_i$, $i=1,2,...,\omega$}.
\end{equation}
The computational complexity of \ref{dehmergraphentropyGT1} and \ref{dehmergraphentropyGT1} is classified as polynomial.

\subsection{Some other entropy measures}
In this subsection, we review some other well-known entropy measures. For discrete-time Markov chains \(X=\{X_t\}\) and \(Y=\{Y_t\}\), $t=1,2,...$, \textbf{Transfer entropy} \citep{tranferentropybook} quantifies the decrease in the amount of uncertainty in future values of stochastic process \(X\) conditioning on the past values of \(X\) itself by knowing the past values of another stochastic process \(Y\). It measures how much additional information the source process \(Y\) provides about state transitions in the target process \(X\) beyond what is already explained by \(X\)'s past. It is defined as 
\begin{equation}
H_{T}^{m,n}(X,Y;t) = H(X_t|X_{t-1}(m))-H(X_t|X_{t-1}(m),Y_{t-1}(n)),
\end{equation}
where $H(X_t|X_{t-1}(m))=H(X_t,X_{t-1}(m))-H(X_{t-1}(m))$ is a conditional Shannon entropy and $X_{t}(m)\equiv (X_{t},X_{t-1},...,X_{t-m+1})$. $H_{T}^{m,n}(X,Y;t)$ represents the uncertainty in current state \(X_t\) that is removed by both the past states of \(X_t\) and \(Y_t\). $H_{T}^{m,n}(X,Y;t)=0$ if and only if the state $X_t$ of \(X\), given its past $X_{t-1}(m)$, is independent of the past states $Y_{t-1}(n)$ of \(Y\). The assumption of a Markov chain is made for simplification. We can also express it in terms of mutual information as 
\begin{equation}
    H_{T}^{m,n}(X,Y;t)=I(X_t;Y_{t-1}(n)|X_{t-1}(m)).
\end{equation}

The spectrum of a signal $X$ describes how its power is spread over different frequencies, illustrating the energy content at each frequency. Computing the power spectral density 
\begin{equation}
P(X=x_i)=\frac{1}{N}{|X(x_i)|}^2
\end{equation}
of \( X \), where \( N \) is the bin number used in the frequency computation. Normalize the power spectral density $P(X=x_i)$ to derive a pmf as 
\begin{equation}
    p_i=\frac{P(X=x_i)}{\sum_j P(X=x_j)}.
\end{equation}
Thus, the \textbf{spectral entropy} \citep{inouye1991quantificationabcdef} of a signal \( X \) is defined as 
\begin{equation}
    H_{SpE}(X)=-\sum_{i}p_i\ln(p_i).
\end{equation}

When the distribution or frequency of an event is uncertain or discrepancies exist in the data, Dempster-Shafer evidence theory \citep{dempsterevidencetheory2008upper} plays an important role in assessing the associated uncertainty of the sample. In such cases, the associated uncertainty is quantified by \textbf{Deng entropy}. Consider \( X \) as a set of events that are mutually exclusive and collectively exhaustive, and \( m \) is the corresponding pmf on the power set of \( X \). The Deng entropy \citep{deng2016dengentropy} of \( X \), based on \( m \), is defined as 
\begin{equation}
\begin{split}
    H_{DgE}^{m}(X)&=-\sum_{Y\subseteq X}m(Y)\log_{2}\left(\frac{m(Y)}{{2^{|Y|}-1}} \right)\\
    &=\sum_{Y\subseteq X}m(Y)\log_{2}\left({{2^{|Y|}-1}} \right)-\sum_{Y\subseteq X}m(Y)\log_{2}\left({m(Y)} \right).
    \end{split}
\end{equation}
Here, the first term corresponds to the non-specificity in the pmf $m$, while the second term is typical for entropy measures. Considering singleton subsets of $X$ only, we have $H_{DgE}^{m}(X)=H(X)$, the Shannon entropy (\ref{shannonentropy}), which implies $H_{DgE}^{m}$ satisfies the probability consistency property\citep{abellan2017analyzingpropertiesofdengentropy}. \( H_{DgE}^{m} \) is not additive because the non-specificity component of the entropy measure is non-additive, while the Shannon entropy component is additive. It is also referred to as \textbf{Belief entropy}.

Let \( X \) be a compact topological space and \( \mathcal{U} \) an open cover of \( X \). Let \( N(\mathcal{U}) \) be the cardinality of minimal subcover of \( X \) contained in $\mathcal{U}$. Let \( \Phi: X \to X \) be a continuous mapping. Define
\begin{equation}
\mathcal{U} \wedge  \mathcal{V} = \{ A\cap B: A\in \mathcal{U}, B\in \mathcal{V} \},
\end{equation}
for any two open covers \( \mathcal{U} \) and \( \mathcal{V} \). The entropy of the cover \( \mathcal{U} \) is given by \( E_{T}(\mathcal{U})=\ln(N(\mathcal{U})) \) and the topological entropy of $\Phi$ with respect to cover $\mathcal{U}$ is defined as 
\begin{equation}
    h(\Phi, \mathcal{U})=\lim_{n \to \infty} \frac{E_{T}(\mathcal{U} \wedge {\Phi^{-1}{\mathcal{U}}} \wedge \cdots \wedge {\Phi^{-n+1}{\mathcal{U}}} )}{n}.
\end{equation}
\textbf{Adler topological entropy} \citep{adler1965adlertopologicalentropy} of the mapping \( \Phi \) is defined as 
\begin{equation}
    H_{ATE}(\Phi)=\sup_{\mathcal{U}}h(\Phi, \mathcal{U}).
\end{equation}
$H_{ATE}$ always takes values on the extended real line. The entropy function $H_{ATE}$ is invariant, meaning 
\begin{equation}
 H_{ATE}(\Psi\Phi{\Psi^{-1}}) = H_{ATE}(\Phi) , 
\end{equation}
where ${\Psi}$ is a homeomorphism from \( X \) onto some \( X^{'} \). For a positive integer $n$, we have 
\begin{equation}
    H_{ATE}(\Phi^{n})=n H_{ATE}(\Phi).
\end{equation}
Also, if $\Phi$ is a homeomorphism, then for any integer $n$, we have 
\begin{equation}
    H_{ATE}(\Phi^{n})=|n| H_{ATE}(\Phi).
\end{equation}
Suppose \( X \) and \( Y \) are two compact topological spaces, and \( 
\Phi_1 \) and \( \Phi_2 \) are two continuous mappings from \( X \) and \( Y \) to themselves, respectively. Additionally, \( \Phi_1 * \Phi_2 \) is a continuous function from \( X * Y \) to itself by 
\begin{equation}
    \Phi_1 * \Phi_2 (x,y)=(\Phi_1(x),\Phi_2(y)).
\end{equation}
Then, we have 
\begin{equation}
H_{ATE}(\Phi_1 * \Phi_2) = H_{ATE}(\Phi_1) + H_{ATE}(\Phi_2).
\end{equation}

Consider the pdf \( g_{\mathcal{X}}(x) \) and \( g_{\mathcal{Y}}(y) \) corresponding to the random variables $\mathcal{X}$ and $\mathcal{Y}$, respectively. The \textbf{cross-entropy} (or Kullback-Leibler divergence) \citep{de2005tutorialcrossentropy} between \( \mathcal{X} \) and \( \mathcal{Y} \) is defined as 
\begin{equation}
    H_{CrE}(X:Y)=\int_{-\infty}^{\infty} \int_{-\infty}^{\infty} g_{\mathcal{X}}(x) \log \left( \frac{g_{\mathcal{X}}(x)}{g_{\mathcal{Y}}(y)} \right)dxdy.
\end{equation}
Cross entropy has key properties that make it useful for comparing probability distributions. It is always non-negative, meaning it is zero only when two distributions are identical. However, it is not symmetric, as the divergence from one distribution to another is not the same when reversed. Furthermore, it is not a true metric because it lacks the properties of distance measures, such as symmetry and triangle inequality.

\section{Applications of Entropy measures}\label{section3}

In this section, we focus on the applications of entropy measures discussed in the previous section with their benefits, and provide a comprehensive view for different applications domains.

\begin{figure}
    \centering
    \includegraphics[width=1\linewidth]{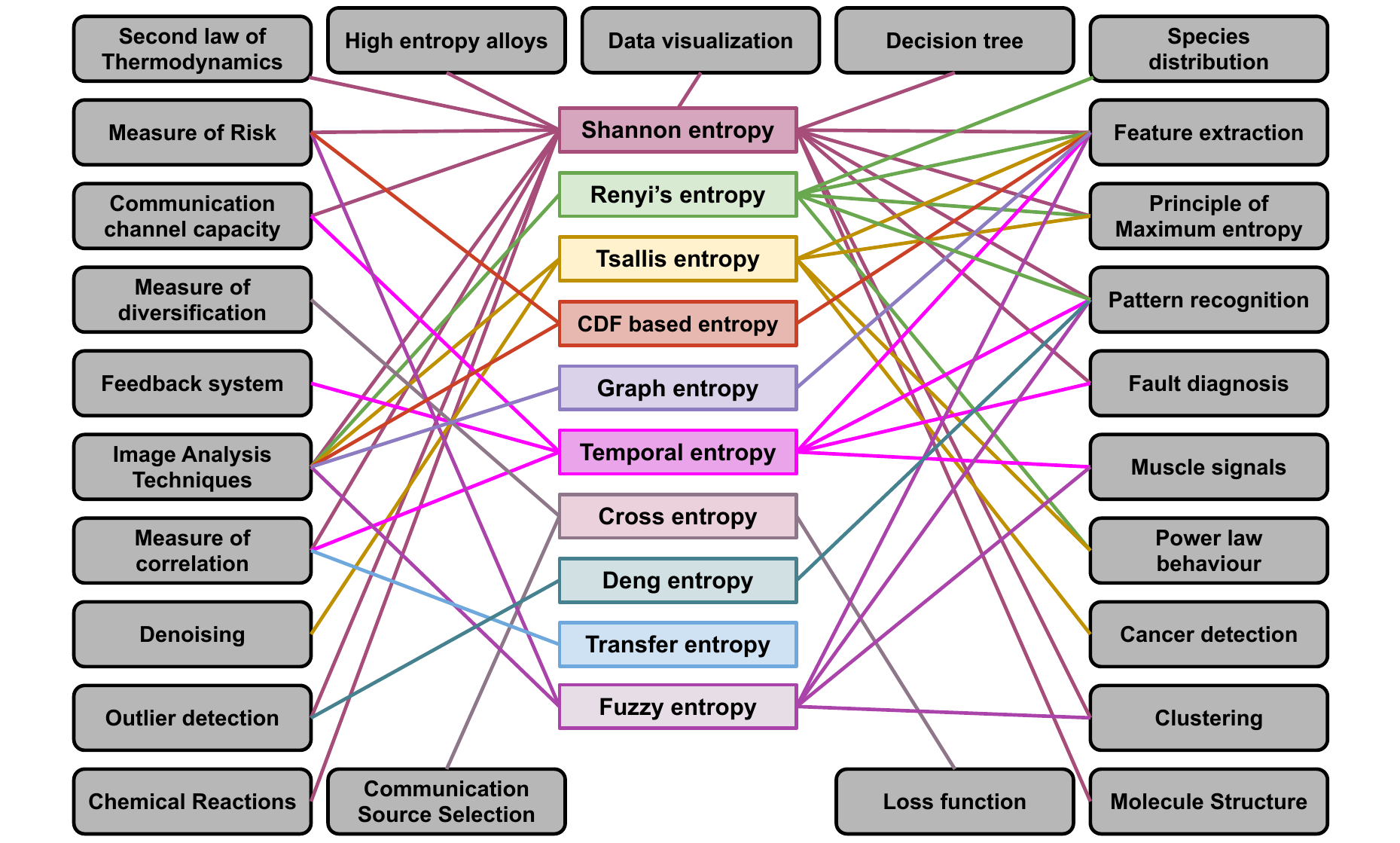}
    \caption{Some applications of entropy measures discussed in this section}
    \label{Thermodynamic systems}
\end{figure}

\subsection{Entropy in Thermodynamics}

Thermodynamics examines the principles of heat, temperature, and energy conversion, exploring how energy is exchanged and transformed in physical systems. Energy \citep{balibrea2016clausiusthermodynamics1} is neither generated nor annihilated but can be transformed into equivalent mechanical work, which means that total energy must remain constant. Clausius demonstrated that the total change in entropy over a complete reversible cycle for an ideal system is zero, meaning the system gains as much entropy when it absorbs heat as it loses when it cools. Boltzmann \citep{boltzmann1974secondboltmannthermodynamics1} rigorously established a precise relationship between gas's temperature and the average kinetic energy of its constituent molecules, illustrating that thermal energy and entropy are correlated with this kinetic movement. In an equilibrium state, where entropy reaches its maximum, and no heat exchange occurs between substances of equal temperature, exchange of kinetic energy stops. The second law of thermodynamics\citep{moran2010fundamentals2ndlawofthermodynamics1} can also be defined as a closed system's entropy never diminishes, regardless of the processes occurring within it: \(\Delta S_{B}^{W} \geq 0\), where \(\Delta S_{B}^{W} = 0\) corresponds to reversible processes, while \(\Delta S_{B}^{W} > 0\) represents irreversible processes. Here \(\Delta S_{B}^{W}=\frac{\Delta Q}{ T} \), where $Q$ is the heat transferred over temperature $T$.

 Wehrl \citep{wehrl1978generalapplicationonthermodynamics1} described Shannon entropy as a quantitative measure of the chaotic properties of a microscopic system. It bridges the macroscopic and microscopic realms of nature, elucidating the behaviour of macroscopic systems—such as real matter—in states of equilibrium or near equilibrium. The third law of thermodynamics posits that for systems with nondegenerate ground states in equilibrium, change in entropy must tend to zero as the temperature tends to absolute zero, i.e., 0 K. Conversely, entropy can be precisely zero exclusively at absolute zero temperature. Bento \citep{bento2015thirdlawthermodynamics1} investigated the third law of thermodynamics in the context of Tsallis and Kaniadakis entropy measures. The study delineated the conditions under which the third law of thermodynamics is valid.

 \begin{figure}
    \centering
    \includegraphics[width=0.9\linewidth]{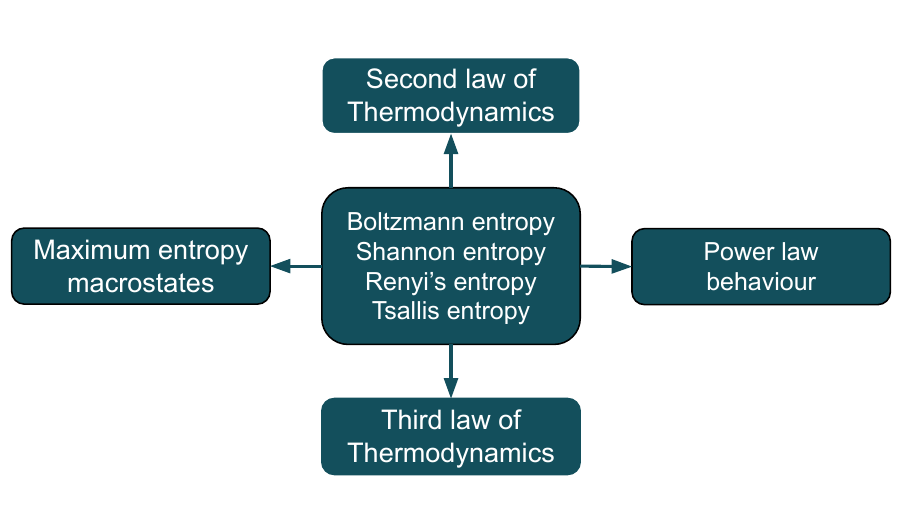}
    \caption{Entropy Measures and their Applications in Thermodynamic Systems}
    \label{Thermodynamic systems}
\end{figure}

 Jaynes \citep{kapur1992spectralentropy} uses Shannon's framework to establish the PME, declaring that the most impartial probability assignment maximizes entropy given the constraints of the available information. Further, the probability distributions such as normal, exponential, gamma and Nakagami are derived uniquely by maximizing Shannon entropy under distinct constraints\citep{kapur1992spectralentropy,kumar2024q}. The PME is related to determining the most probable macrostate of a system, which is a foundational concept of statistical mechanics. Through the PME \citep{bashkirov2004maximumMrenyithermodynamics1, martinez2000Mtsallisthermodynamics1, macedo2013maximumMkaniadakisthermodynamics1} with measures such as Shannon, Rényi, Tsallis, and Kaniadakis, this framework enables the computation and generalization of state distributions of a system, capturing long and short-range correlations and allowing the modelling of power-law behaviour in systems.

 The concept of entropy is used in constructing high entropy alloys\citep{george2019definationofhighentropyalloysthermodynamics1}, a newly discovered research area in material science. Typically, these are formed by using five or more different atoms in equal proportion\citep{miracle2017criticalhighentropyalloythermodynamics1}. The lattice position of the atoms estimates the configurational entropy. The atom's positions generate a high entropy, which considerably improves the alloy's microstructure. These alloys exhibit features such as high strength, flexibility, corrosion resistance, and thermal stability. This makes them valuable in the aerospace and automotive industries.

\subsection{Entropy in Communication Theory}
In communication theory\citep{continuousshannonentropycover1999elements}, the main objective is to send a message through a communication channel. The sender transmits a series of partial messages that provide hints about the original message. The information content of each partial message indicates how much uncertainty it reduces for the receiver. In this context, entropy represents the average number of bits needed to describe each message, taking into account all possible messages that can be sent. Shannon\citep{shannon2001mathematical} laid the foundation of information theory, and to fully understand applications of his results in communication theory, we need an understanding of the concepts of the source, channel, and transducer. The origin of the information that will be transmitted is known as the source. The channel facilitates the transmission of this message from the source to the receiver. A device transducer converts energy or signals from one form to another. It is considered non-singular if each input results in a unique output, ensuring no loss of information during the conversion process.

\begin{figure}
    \centering
    \includegraphics[width=0.6\linewidth]{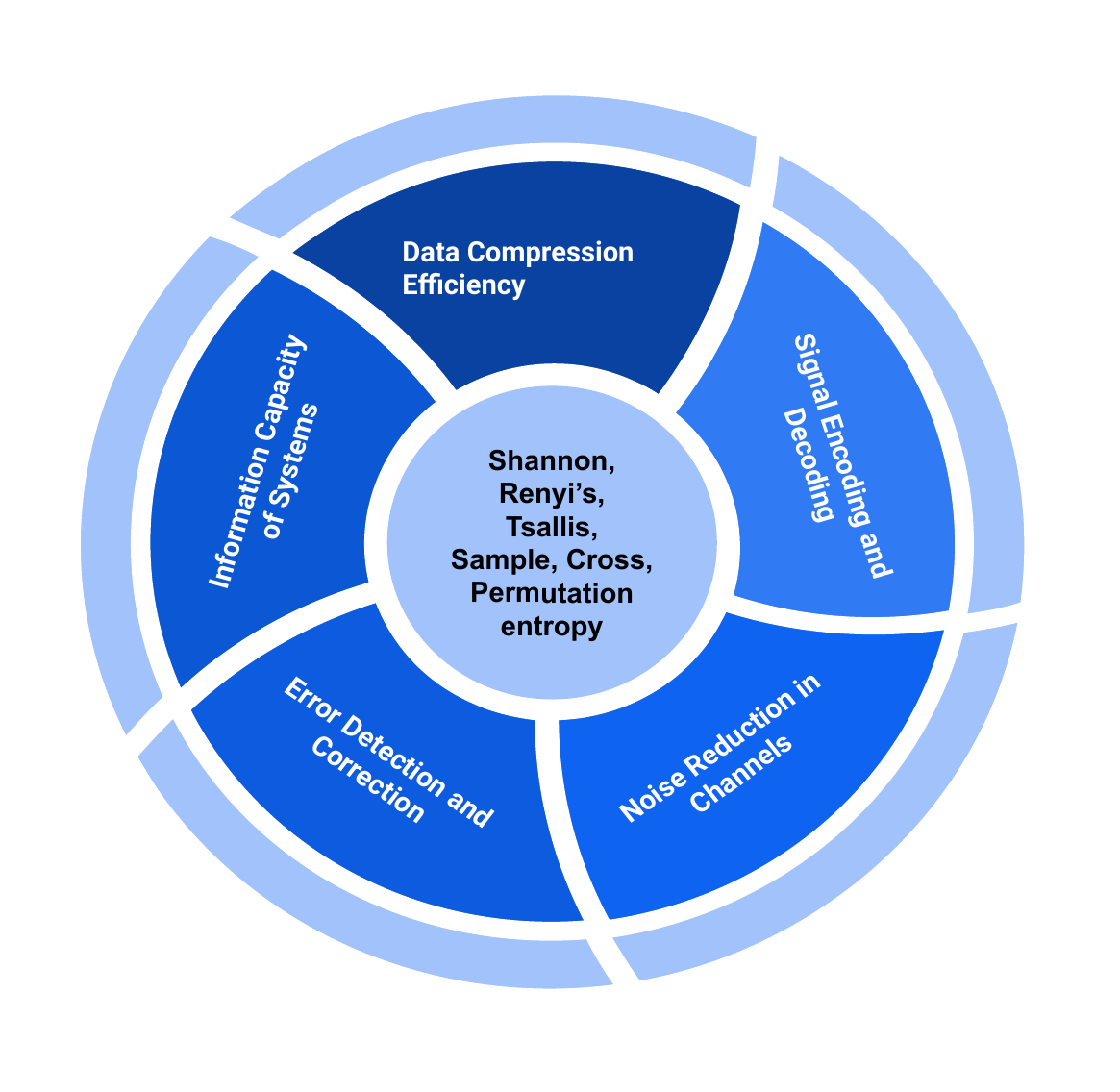}
    \caption{Entropy Measures and their Applications in Communication Theory}
    \label{Communication theory}
\end{figure}

\begin{enumerate}
    \item On processing input from a finite state statistical source, the output of a finite state transducer is a finite state statistical source with the output's entropy equals or is less than the input's (per unit time). The input and output entropies are the same if the transducer is non-singular.
    \item If a source has an entropy \( \mathcal{H} \) (bits per symbol) and a channel has a capacity \( \mathcal{C} \) (bits per second), then it is possible to encode the source's output so that the message is transmitted at an average rate of \( {\mathcal{C}}/{\mathcal{H}} \) symbols per second, with only a very small error \( \varepsilon \), and it is not possible to transmit at a rate higher than \( {\mathcal{C}}/{\mathcal{H}} \).
\end{enumerate}

Rényi's entropy\citep{csiszar1995generalizedRenyisentropycommunication2} refines the analysis of communication channels by extending Shannon's entropy to account for varying error rates and coding structures. It captures mutual information across different orders, offering a more nuanced approach for complex coding channels and hypothesis testing in constant composition codes. Tsallis entropy\citep{mukherjee2019performanceTsalliscommunication2} describes heavy tail in log-normally distributed data with a non-extensive parameter $q$ through the q-log normal distribution, and for $q = 1.8$ for a simulated dataset, it improves fading, outages, and channel capacity. SE\citep{holliday2006capacitySampleentropycommunication2} and permutation entropy\citep{xiang2014phasepermutationentropycommunication2,wang2018securitypermutationentropycommunication2} are used to evaluate chaotic signals and characterize channel capacity, offering improved insights into feedback systems using Lyapunov exponents and auto-correlation analysis, respectively. In cloud computing systems, cross-entropy\citep{rahmani2020kullbackleiblerdivergencecommunication2} reduces migration and improves energy efficiency in selecting source and destination hosts.

\subsection{Entropy in Financial Markets}
The study of finance focuses on managing money, assets, and other financial means to optimize wealth creation, risk control, and resource distribution. Financial studies focus on topics such as market behaviour, corporate finance, and economic stability. Philippatos and Wilson\citep{philippatos1972entropyshannonfinance3} were the first to apply Shannon entropy in portfolio construction. They used a mean-entropy approach, comparing it with the Markowitz and Sharpe models, and found that the results were consistent across the methods. A fuzzy mean entropy model, utilizing fuzzy entropy calculated via the credibility function, is presented in \citep{huang2008meanfuzzyentropyfinance3}. This model maximizes the expected mean of a portfolio while accounting for fuzzy entropy uncertainty. A comparative analysis with the fuzzy mean variance model demonstrates its effectiveness. 

\begin{figure}
    \centering
    \includegraphics[width=0.62\linewidth]{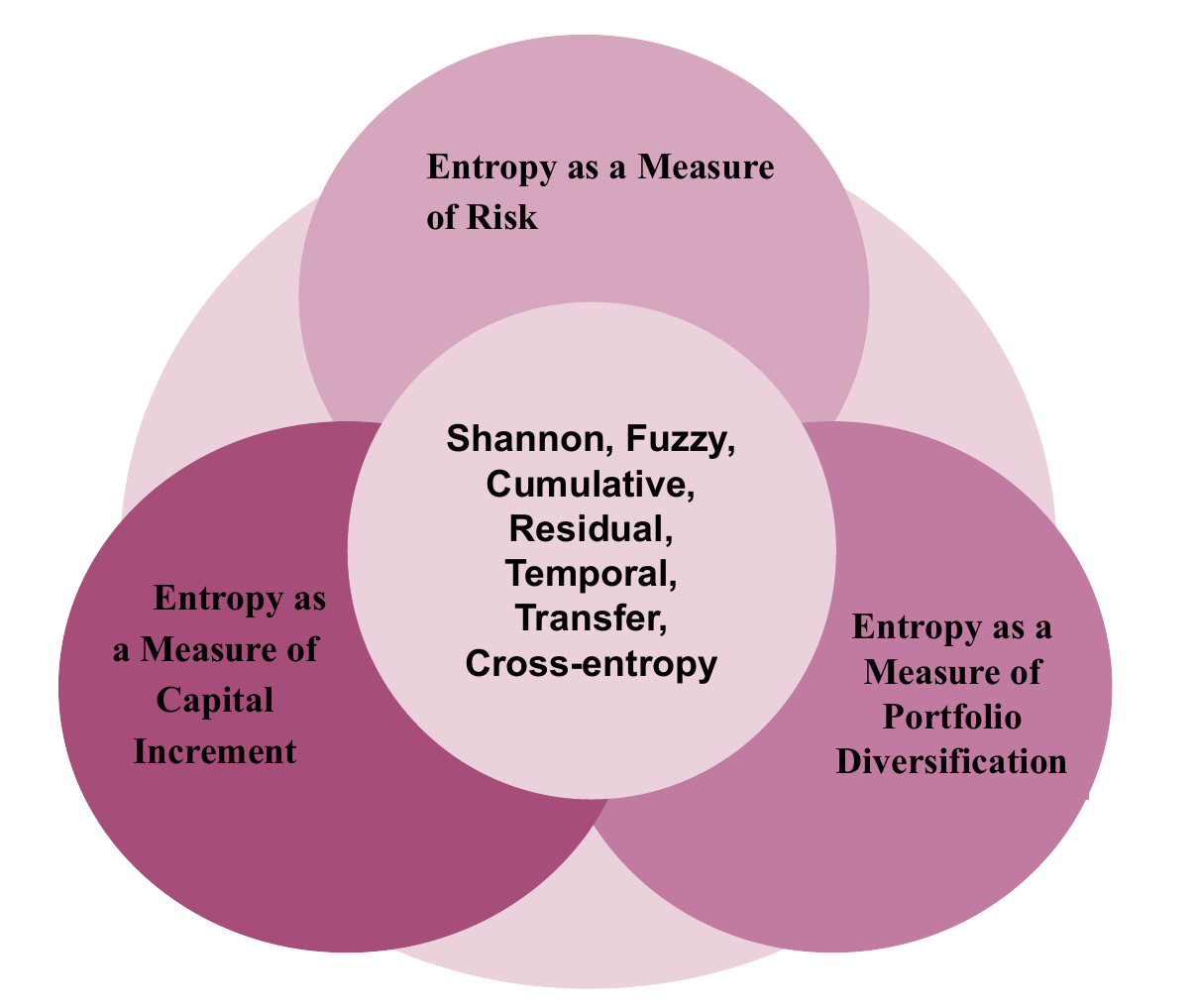}
    \caption{Entropy Measures and their Applications in Financial Markets}
    \label{Financial markets}
\end{figure}

A CE\citep{baghmolaei2022applicationCRE1finance3} approach has been applied in place of variance within the Markowitz model, leading to a reduction in risk associated with portfolio selection. CRE further serves as an effective tool for quantifying uncertainty in stocks, particularly for heavy-tailed distributions\citep{yang2012studyCRE2finance3}, providing reliable descriptions even in cases, where variance is undefined. A comparative analysis of various entropy measures\citep{zhou2017propertiesCRE3finance3}, including Shannon, fuzzy, and cumulative residual entropy, reveals that the fuzzy entropy model outperforms others in terms of daily and relative cumulative returns. Additionally, a variant of CRE has been developed to study asset risk\citep{sun2022statisticalCRE4finance3} under extreme market conditions, evaluating the influence of a broader set of stocks.

Multiscale entropy quantifies the complexity of financial signals, categorizing them based on non-linear correlations measured through this approach. Variants of permutation entropy are more adept at accurately distinguishing multiple scales\citep{xia2012multiscaleentropyfinance3,yin2014weightedpermutationentropyfinance3}. Among these measures, approximate, sample, and dispersion entropy effectively categorize stocks by the complexity of their financial signals, as detailed in \citep{olbrys2022approximatesampleentropyfinance3} and \citep{wang2021generalizeddispersionentropyfinance3}. To meet specified conditions on asset weights and associated risk, a cross-entropy model\citep{bera2008optimalcrossentropyfinance3} is proposed that uses cross-entropy as the objective function, with constraints based on the mean and variance-covariance matrix. The reference distribution in the cross-entropy formulation is chosen on the basis of the desired conditions, such as achieving an equally weighted portfolio or a minimum variance. Transfer entropy\citep{dimpfl2013usingTransferentropyfinance3}, a model-free measure, captures the information flow between different stocks using transition probabilities without being constrained by linear dynamics.

We have not included numerous other models based on entropy, which are used to quantify risk, diversify portfolios, and manage the uncertainty associated with risky assets, due to space constraints.

\subsection{Entropy in Categorical Data Analysis}
Data from sources influenced by multiple factors are often categorized to highlight the impact of each factor. In an $n$-dimensional table, $n$ categorical variables influence the outcome, with a scale to distinguish between different categories and their effects. 

Outlier detection in categorical data is a crucial task, particularly in high-dimensional datasets. Wu\citep{wu2011informationShannonentropycategoricaldata1} introduced two one-parameter algorithm that utilizes Shannon entropy and correlation to identify outliers, assigning likelihoods to each detected outlier. These methods perform effectively on large-scale and high-dimensional data, with results demonstrated on both real and simulated datasets. Visualizing categorical data is difficult because of its discrete nature. However, by utilizing the Shannon entropy of marginal and joint categorical variables, the study\citep{alsakran2014usingshannonentropycategorical2} provides an effective approach to categorical data visualization. It efficiently manages high-dimensional data, delineates its boundaries, and allows for testing and tuning of variables to enhance the visualization process. Based on Shannon entropy, a monitoring technique\citep{das2017detectingshannoncategorical3} is developed to quantify uncertainty in contingency tables using non-parametric estimation and a dependency measure\citep{skotarczak2018entropyshannoncategorical4} for categorical variables. 

In categorical tables, a key task is the grouping of data points with similar patterns, properties, or characteristics, where clustering techniques offer efficient and accurate decision-making. A modified form of Shannon entropy is introduced as a measure\citep{simovici2000generalizedshannoncategorical5} for clustering in the categorical table, along with their impurity. This approach offers flexibility and, in some cases, serves as a natural average distance for clustering. COOLCAT\citep{barbara2002coolcatcategorical6} is a Shannon entropy-based categorical data clustering algorithm that minimizes the Shannon entropy within clusters, outperforming algorithms like ROCK. It demonstrates stability across various domain samples, scales effectively to large datasets, and is evaluated using a categorical utility function. In both hard and fuzzy clustering algorithms\citep{mahamadou2020categorical7}, where hard clustering relies on certainty for cluster inclusion and fuzzy clustering uses a membership function, Shannon entropy is integrated into the objective function. Entropy quantifies the weights corresponding to the clusters during the minimization of the objective function. Fuzzy entropy clustering algorithm\citep{mahamadou2020categorical7} was tested on ten real-life categorical datasets, yielding favourable results.

\begin{figure}
    \centering
    \includegraphics[width=0.6\linewidth]{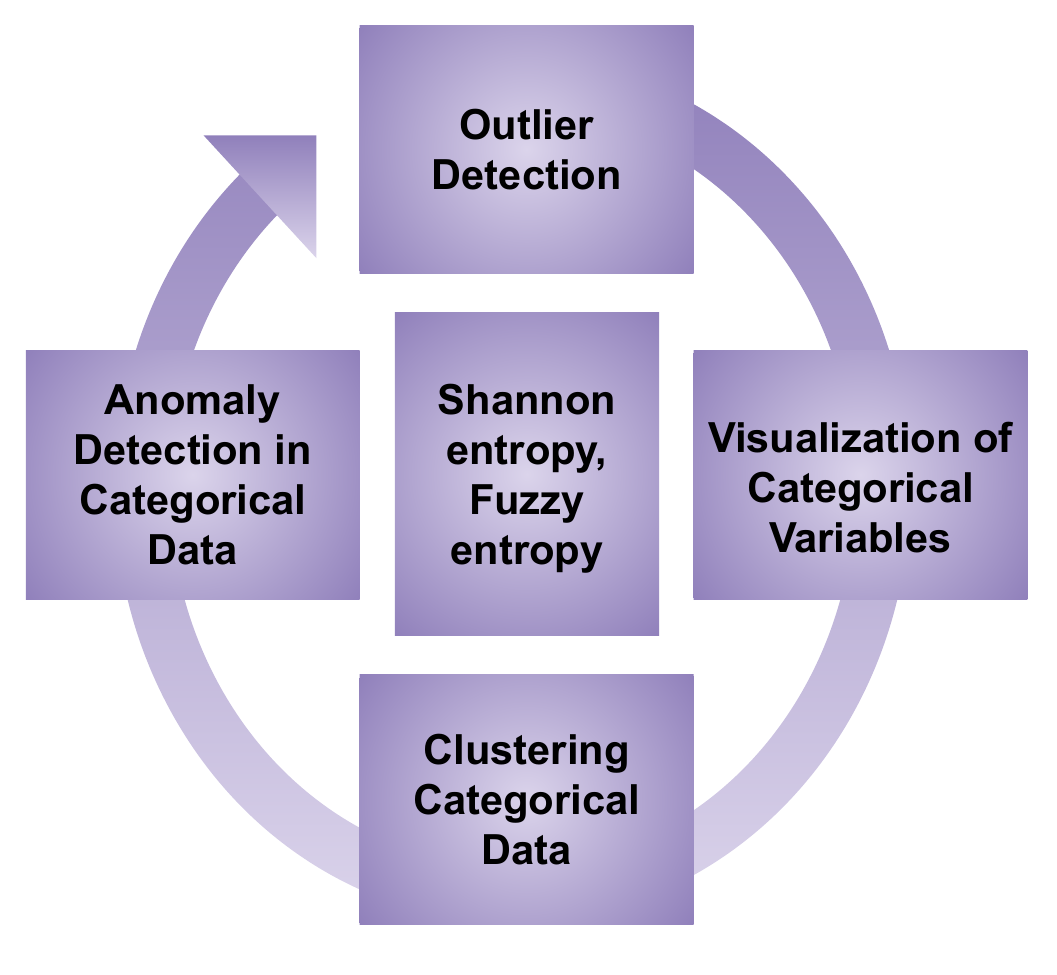}
    \caption{Entropy Measures and their applications in Categorical data}
    \label{Categorical data}
\end{figure}

Anomaly detection in categorical data is often based on identifying responsible patterns, with detection involving the computation of distance from these patterns. Points at significant distances are labelled as anomalies. Techniques addressing such anomalies in categorical data utilize Shannon entropy and correlation-based measures, as detailed in \citep{taha2019anomalycategorical8}.

\subsection{Entropy in Artificial Intelligence}
This sub-section is organized into two distinct parts: the first part focuses on the entropy application in image processing, while the second part explores its utilization in various machine-learning techniques.

\subsubsection{Entropy in Image Processing}
Image processing is vital in research as it enables the extraction, analysis, and augmentation of visual data needed for precise interpretation and decision-making. It is broadly utilized in medical diagnostics, remote sensing, and material research to uncover concealed patterns, improve picture quality, and enable data-driven insights. We outline some applications of entropy in image processing.

\begin{figure}
    \centering
    \includegraphics[width=0.7\linewidth]{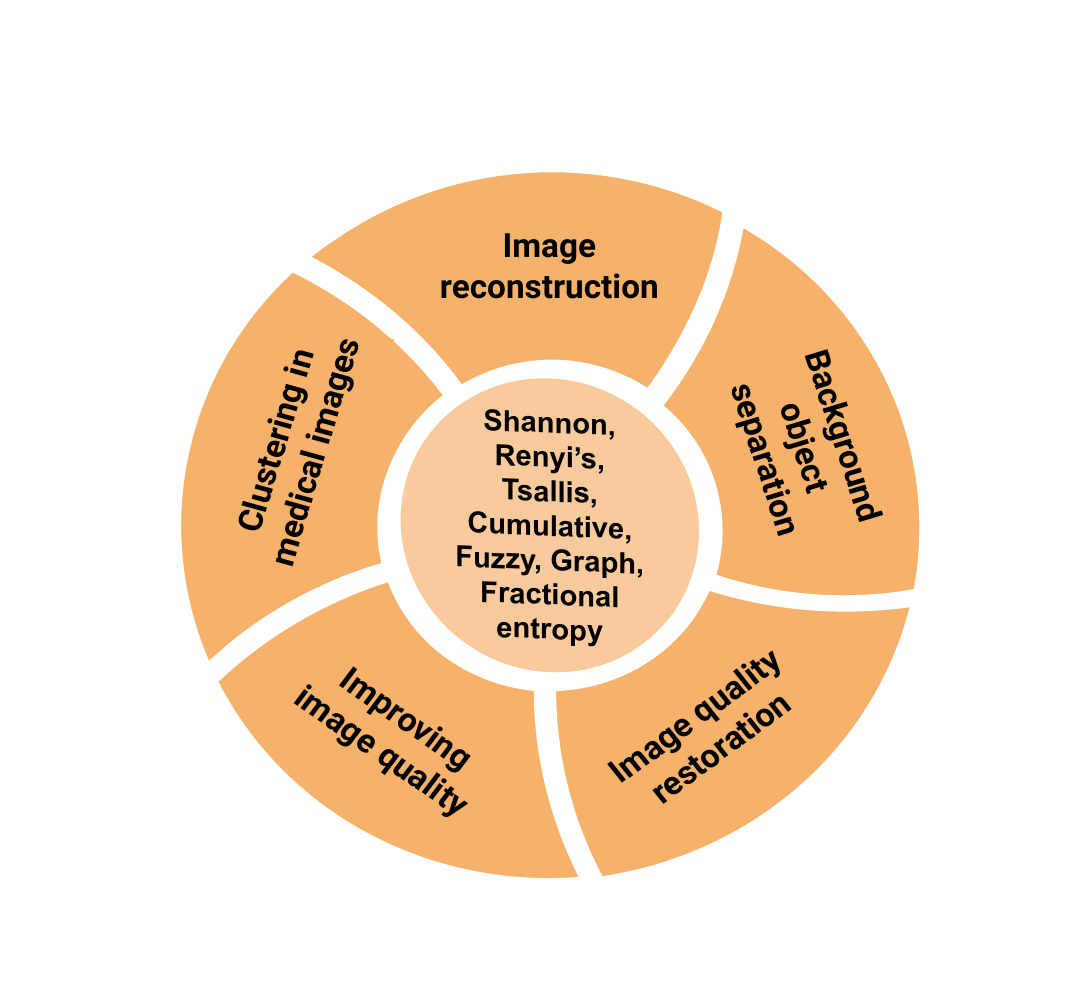}
    \caption{Entropy Measures and their Applications in Image Processing}
    \label{Image processing}
\end{figure}

\begin{enumerate}
    \item Shannon entropy, is used to reconstruct images from incomplete or noisy data by incorporating prior knowledge about the object\citep{gull1984maximumIMG11}, improving image quality. It effectively detects edges by measuring brightness and colour entropy in small areas\citep{shiozaki1986edgeIMG12}, identifying sharp changes with smooth edge detection in grayscale and colour images. Measuring it locally also improves randomness detection in shuffled and encrypted images, offering a more accurate and insightful analysis than the traditional global method\citep{wu2013localIMG13}.

   \item Generalized entropy measures, such as Renyi's, Tsallis, and Kaniadakis entropy, extend Shannon entropy's principles to improve image processing tasks like segmentation. These methods\citep{sahoo2004thresholdingIMG21,de2004imageIMG22,sparavigna2015roleIMG23,lei2020adaptiveIMG24} are used to separate objects from the background by incorporating additional information, such as local pixel values and non-additive properties, which Shannon entropy alone may not account for. Each form of entropy offers specific advantages for different types of images, including medical images and those with long-tailed distributions, improving accuracy and performance in identifying objects within images.

   \item CRE is applied in multi-modality image registration\citep{wang2007nonIMG31}, offering faster and more reliable alignment of images with different contrasts and brightness levels compared to other methods. Additionally, it enhances the fusion of high-resolution colour images with lower-resolution hyperspectral data\citep{hasan2010regisrationIMG32}, demonstrating its utility in surveillance image alignment. In image denoising, an improved auto-encoder model based on CRE and residual statistics is introduced\citep{xiang2020imageIMG33}, achieving better noise suppression and image quality restoration. A new edge detection method using cumulative residual entropy is also proposed\citep{al2020novelIMG34}, showing superior performance over traditional techniques, as evidenced by improved peak signal-to-noise ratios across various image types.

   \item As mentioned earlier, fuzzy entropy is a measure used to assess the sharpness of images, where its value increases as the image becomes blurrier. This concept improves image thresholding by defining a ``sharper than" relationship between fuzzy sets, facilitating better separation of objects from backgrounds in degraded images. Additionally, the authors introduce a genetic algorithm\citep{di1998imageIMG41,sanyal2011adaptiveIMG42,yin2017unsupervisedIMG43} to automatically select the optimal fuzzy regions for membership functions, enhancing image quality. The authors develop various advanced algorithms\citep{naidu2018shannonIMG44}, such as adaptive bacterial foraging and a firefly algorithm, to optimize fuzzy entropy for more effective image segmentation, demonstrating improved results for both grayscale and colour images compared to traditional methods.

   \item A graph entropy-based method\citep{zhan2016graphIMG51,zhan2015medicalIMG52} for clustering medical brain images allows doctors to identify similar pathology images more efficiently, which aids in disease analysis. Additionally, the introduction of SampEn2D, a SE of the two-dimensional method\citep{silva2016twoIMG53}, demonstrates its effectiveness in distinguishing different textures and accurately classifying biological images, showing its potential as a reliable tool in biomedical image analysis. Furthermore, the development of multiscale entropy for one and two dimensions extends the concept of multiscale entropy to images\citep{silva2018twoIMG54}, effectively analyzing irregularities and performing well across various applications. For example, the Fractional Entropy model\citep{al2018newIMG55} enhances kidney images by identifying edges and improving the quality of MRI images, outperforming traditional methods.
\end{enumerate}

\subsubsection{Entropy in Machine Learning}

Attribute selection\citep{wang1984analysisdecisiontreeAI2} for decision trees is one of the well-known applications of Shannon entropy, where entropy is used as an information measure within the data. Researchers have also conducted a comparative study between Shannon, Rényi, and Tsallis entropy to determine the most efficient tree construction methods\citep{lima2010decisiontreeAI1}. Various learning techniques, such as feature extraction\citep{guido2018tutorialfeature1AI1,sharma2016feature2AI1}, pattern recognition\citep{watanabe1981patternAI1}, sequence complexity analysis\citep{romero2001sequencecomplexityAI1}, group diversity assessment\citep{balch2000hierarchicgroupdiversityAI1}, sentence representation\citep{arroyo2019unsupervisedsentencerepresentationAI1}, fault diagnosis\citep{kankar2011rollingfaultdiagnosis1AI1,li2018entropyfaultdiagnosis2AI1}, and signal classification\citep{sabeti2009entropycomplexitymeasureAI1}, widely utilize Shannon entropy. These applications represent some of the most common uses centred around Shannon entropy.

\begin{table}[h!]
\centering
\resizebox{0.85\textwidth}{!}{%
\begin{tabular}{|c|c|c|c|}
\hline
\textbf{Entropy} & \textbf{Accuracy} & \textbf{False Positive} & \textbf{False Negative} \\ \hline
  Shannon entropy    & $99.638\%$      & $0.8932\%$      & $1.4462\%$      \\ \hline
Renyi's entropy($\alpha=0.5$)     & $99.691\%$      & $0.7656\%$      & $0.7231\%$      \\ \hline
Tsallis entropy($q=1.9$)     & $99.6998\%$      & $0.5104\%$      & $0.6380\%$      \\ \hline
\end{tabular}
}
\caption{Comparative results of Shannon, Rényi, and Tsallis entropy in a decision tree algorithm applied to a traffic dataset, as presented in \citep{lima2010decisiontreeAI1}.}
\label{tab:sample}
\end{table}

In some models, generalized entropies provide greater flexibility and reduce model development processing time. Rényi's entropy\citep{polewski2016combiningGER1AI1} helps in improving active and semi-supervised learning by making training set creation easier and more efficient. When used with the shuffled frog-leaping algorithm, it also enhances crack detection in bridge infrastructure by improving the identification of boundaries and features\citep{abdelkader2021analyzingGER2AI1}. Additionally, Rényi's entropy is a key to information semifields\citep{valverde2021caseGER3AI1}, which aims to mimic how the brain works and help develop stronger artificial intelligence systems. Tsallis entropy aids in computer-aided diagnosis by effectively differentiating pathological brains from healthy ones in MRI scans\citep{zhang2015pathologicalTsallis2AI1}, as well as improving the recognition of isolated objects in image processing\citep{zhang2015pathologicalTsallis2AI1} through advanced feature extraction.

Time-based entropy measures, such as permutation, multiscale, sample, and approximate entropy, assist in modelling fault detection and machine health. These methods are low-cost and effective for capturing patterns in machine behaviour. Since entropy values quantify uncertainty, they are often paired with learning techniques for early fault detection. A comprehensive range of these methodologies is discussed in \citep{huo2020entropyreviewoftemporalentropyAI1}. Temporal entropy, combined with machine learning techniques, efficiently identifies relevant features in signal classification, with multiscale entropy particularly helpful in classifying wireless signals\citep{llenas2017performanceT1AI1} and analyzing human behaviour\citep{fraiwan2021gaugingT2AI1}, for example, enjoyment and visual interest during museum visits.

Fuzzy entropy is used to transform training data, increasing its clarity and enabling users to better understand and work with it\citep{morente2016improvingfuzzy1AI1}. Additionally, it generates unbalanced linguistic label sets, focusing labels where most data points are concentrated, thereby improving classification accuracy. By effectively choosing pertinent characteristics from big, redundant datasets, lowering dimensionality, and filtering noisy data, fuzzy entropy significantly improves recommender systems. Even under data sparsity, this improves recommendation quality and provides accurate and high-speed user predictions. Including fuzzy entropy in deep learning models\citep{saravanan2019fuzzy2AI1} helps to overcome the restrictions of traditional methods and improves performance even further.

Graph entropy is useful for classifying structured data because it captures deep information through subgraph representations\citep{xu2021deepgraph1AI1}, improving the efficiency and accuracy of graph kernels. Its low computational complexity allows it to scale well with large graphs, making it effective for handling more extensive datasets. Graph entropy helps to improve node embedding dimensions in graph representation learning because it includes both structure and feature entropy. Feature entropy connects node features to graph topology, and structure entropy uses normalized degrees to capture higher-order graph structures. Compared to traditional methods, this combination\citep{luo2021graph2AI1} in the Minimum Graph Entropy algorithm significantly improves model performance while reducing computational complexity.

\begin{figure}
    \centering
    \includegraphics[width=0.65\linewidth]{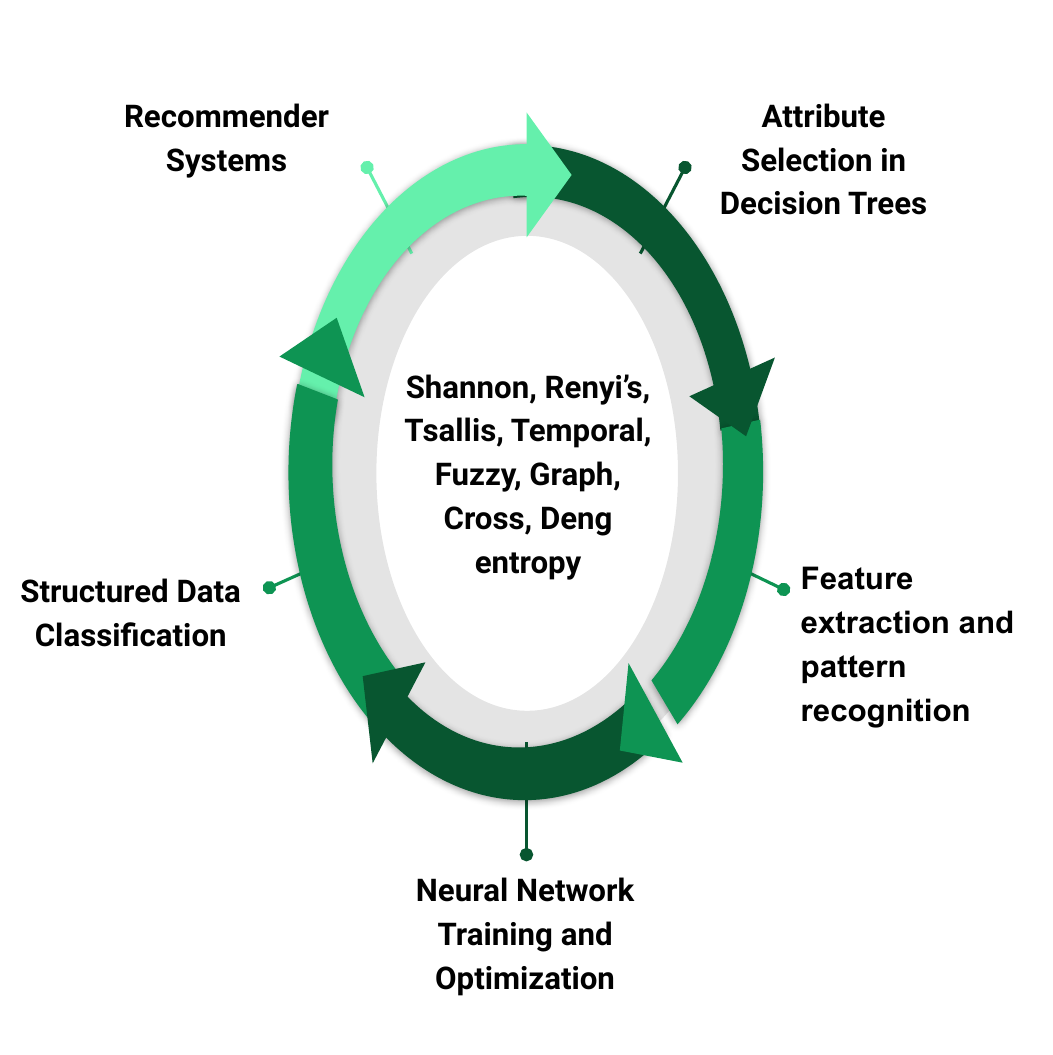}
    \caption{Entropy Measures and their Applications in Artificial Intelligence}
    \label{Artificial Intelligence}
\end{figure}

Cross-entropy is a frequently utilized loss function in machine learning\citep{mao2023cross1AI1}, especially in neural networks to classify tasks, as it aligns with logistic loss when softmax is applied to the outputs. Theoretical analysis shows that cross-entropy belongs to a broader family of loss functions, known as comp-sum losses, which offer robust performance in various settings. A modified cross-entropy loss\citep{ho2019realcross2AI1} improves result by considering real-world costs, such as financial impacts, making it more effective than just focusing on accuracy or F1 score. In neural networks, cross-entropy accelerates the backpropagation algorithm\citep{nasr2002cross3AI1}, improving network performance and reducing training time. Deng entropy is useful for managing complex uncertainty in probability assignments, especially when dealing with more intricate propositions. A modified belief entropy function has been introduced to address Deng entropy's limitations, offering improved accuracy in measuring uncertainty, as demonstrated in pattern recognition examples\citep{cui2019improveddengAI1}.

\subsection{Entropy in Signal Processing Analysis}
Signal processing techniques enable the extraction of valuable patterns from temporal data, while entropy plays a crucial role in quantifying uncertainty and revealing underlying system dynamics, improving the analysis of complex and irregular behaviours. Many applications of Shannon entropy for signal data exist, including source separation, blind deconvolution of autoregressive systems, and density changes (see \citep{bercher2000estimatingshannon1Time1}). Rényi's entropy quantifies signal complexity and distinguishes meaningful neurophysiological activities from noise in time-dependent neuroimaging data\citep{gonzalez2000measuringrenyi1Time1} like fMRI, EEG, and ERP. Calculating Rényi's entropy over time-frequency representations provides a measure of disorder and approximates the number of elementary components in the signal. This approach helps capturing the intricate details of brain activity, offering a deeper understanding of signal structure. 

Tsallis entropy is also used to effectively remove noise from seismic engineering seismograms\citep{beenamol2012waveletTsallis2Time1}, and enhance contrast in mammograms for early cancer detection. In critical applications, a Tsallis entropy-based fuzzy algorithm\citep{kalra2010automaticTsallis1Time1} improves signal clarity and detection accuracy. It provides higher efficacy over conventional methods by significantly improving the signal-to-noise ratio and refining the detection of waves, particularly in high-frequency seismograms. When compared to other thresholding methods, it achieved a significantly higher signal-to-noise ratio. This shows how well it can reduce noise and accurately identify waves, making it an important tool for processing complex seismic data.

\begin{table}[h!]
\centering
\renewcommand{\arraystretch}{1.5} 
\resizebox{\textwidth}{!}{%
\begin{tabular}{|c|c|c|c|c|}
\hline
\textbf{Classifier} & \textbf{Sample entropy} & \textbf{Approximate entropy} & \textbf{Fuzzy entropy} & \textbf{Cumulative residual entropy} \\ \hline
Linear discriminant analysis      & $47.83 \pm 7.77$  &  $ 31.61 \pm 7.46$   & $75.00 \pm 11.26$  &   $78.17 \pm 10.58$  \\ \hline

 Extreme learning machine      & $46.89 \pm 5.87$      & $36.45 \pm 6.14$    &  $77.45 \pm 10.37$  &  $83.72 \pm 7.47$  \\ \hline

Support vector machine      & $49.39 \pm 6.23$ &  $36.11 \pm 7.31$  & 
 $ 79.11 \pm 9.18$  &  $83.95 \pm 6.88$  \\ \hline

\end{tabular}
}
\caption{Comparative classification accuracy results based on Sample, Approximate, Fuzzy, and Cumulative Residual Entropy for multifunctional prosthetic devices, utilizing surface electromyography signal data from channels $1$ and $2$, as detailed in \citep{shi2013semgresidual1Time1}.}
\label{tab:sample_table}
\end{table}

Surface electromyography signals are known for their nonlinear and occasionally chaotic behaviour, making nonlinear time series analysis suitable for feature extraction. CRE is employed to capture the features of surface electromyography data\citep{shi2013semgresidual1Time1}, offering lower computational complexity than fuzzy entropy, SE, and AE. Further, a combined approach\citep{zhou2021generalizedresidual2Time1} using permutation and Rényi's entropy based on CRE is proposed for effectively distinguishing stock markets with varying characteristics when applied. Surface electromyography is important for things like rehabilitation and controlling prosthetics because it helps to assess muscle activity without needing invasive methods. Fuzzy entropy\citep{chen2007characterizationfuzzy1Time1} is better than other methods like SE or AE at spotting changes in muscle activity over time, giving a clearer and more accurate understanding of how muscle signals change. This makes it very useful for identifying different muscle movements. Fuzzy entropy, when used with empirical mode decomposition, in the Inherent FuzzyEn algorithm, helps to better analyze signals by identifying overlapping patterns more clearly\citep{cao2017inherentfuzzy2Time1}. The study shows that fuzzy-based approaches work better than other methods like sample or AE, giving more accurate results in real-world signal analysis.

The most commonly used entropies for time-domain data include permutation, approximate, sample, and dispersion entropies. We include some recent articles on the applicability of these measures, highlighting recent applications such as the use of permutation entropy in analyzing economic markets\citep{zanin2012permutationTime1}, AE in biosignal analysis\citep{fusheng2001approximateTime1}, SE for physiological signals\citep{alcaraz2010reviewsampleTime1}, and dispersion entropy in rotary machines\citep{rostaghi2019applicationdispersionTime1}.

\begin{figure}
    \centering
    \includegraphics[width=0.6\linewidth]{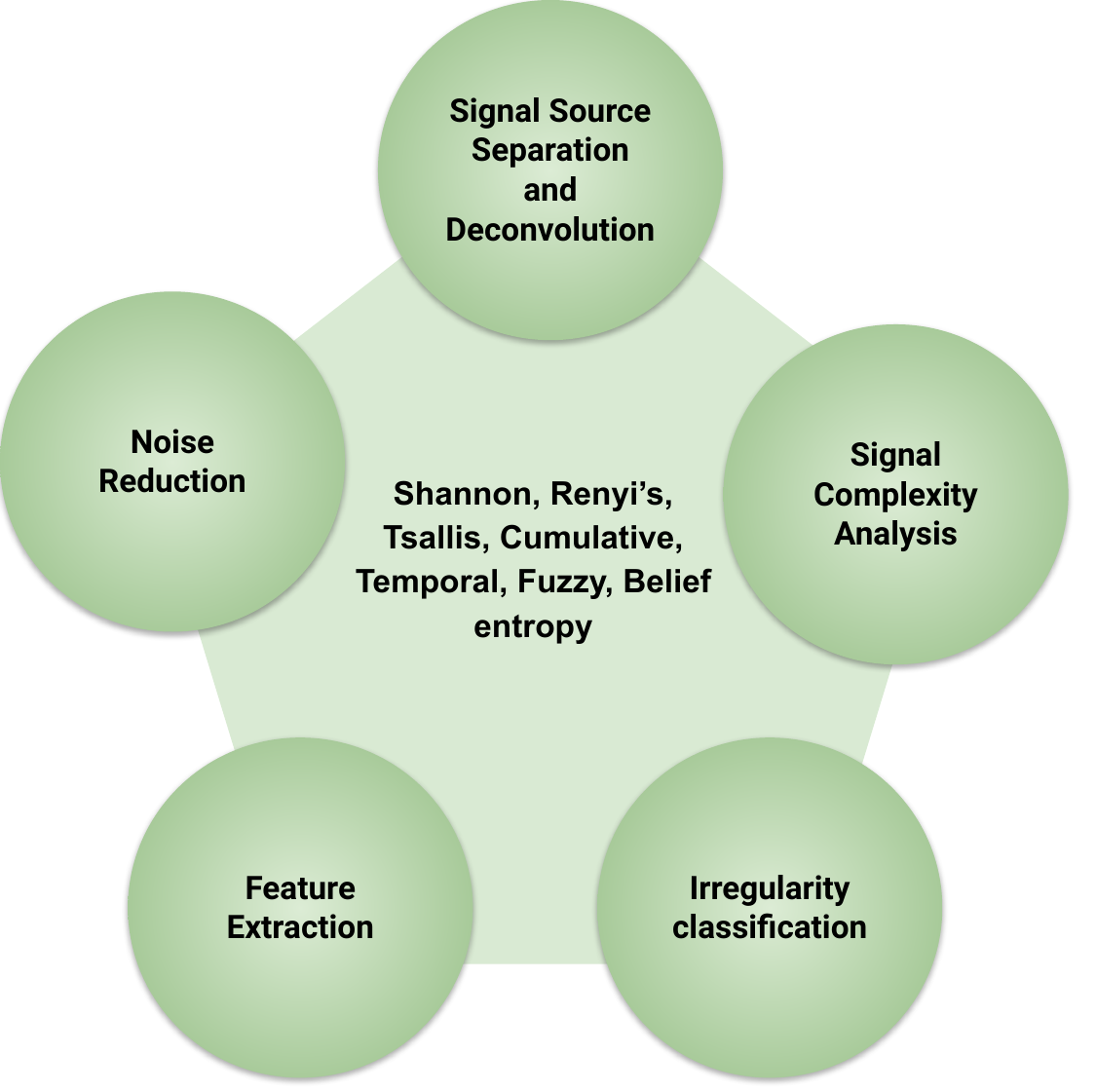}
    \caption{Entropy Measures and their Applications in Time Series and Signal Processing}
    \label{Signal processing}
\end{figure}

Song\citep{song2022combininggraph1Time1} presented a belief entropy-based method for visibility graphs that combines belief entropy with weighted visibility graphs to enhance time series analysis, such as EEG classification. By eliminating outliers and accurately fusing sequence data, it improves classification accuracy, leveraging belief entropy's ability to evaluate internal conflicts in data. Permutation entropy\citep{fabila2022permutationgraph2Time1}, serves as a decisive nonlinear measure of irregularity in time series data, allowing for the quantification of complexity in a variety of forms. By generalizing permutation entropy to analyze signals on irregular graphs, John\citep{fabila2022permutationgraph2Time1} extends its applicability beyond traditional time series and images, preserving its essential properties while enabling new insights into complex data structures. The dispersion entropy approach\citep{fabila2023dispersiongraph3Time1} for graph signals effectively captures intricate signal dynamics, making it valuable for theoretical research and practical applications in multivariate time series and images, including finance, biology, industrial processes, and international events.

\subsection{Entropy in Chemical Processes}
Entropy is important in chemistry for two main reasons: studying molecular graphs and analyzing electron density in molecules. In molecular graphs, entropy helps create models that explain how the structure of a molecule relates to its activity and properties, useful in fields like organic chemistry and drug design. It also helps researchers understand how electrons are distributed in molecules and how this changes during chemical reactions. Additionally, entropy is applied to molecular movements when molecules act as signal carriers. Overall, entropy connects chemistry to other fields like thermodynamics and computer science, opening up new areas of research. In this direction, references \citep{bonchev1995kolmogorovshannontopologyentropychemistry1}, \citep{dehmer2011historygraphentropychemistry1}, and \citep{barigye2014trendsinformationindiceschemistry1} provide a detailed introduction to the applications of Shannon, topological, and graph entropy.

\subsection{Entropy in Biological Systems}
Entropy in biology is used to understand various patterns and processes. Among the most common uses of Shannon's entropy measuring the diversity of living things in an ecosystem, such as the variety of species or types of cells is one. It helps scientists see how different organisms are connected and how they are distributed in space. In studying evolution, entropy explains how systems change over time and how complex populations are organized. Additionally, it can measure how efficiently living things use energy in their metabolism. Overall, entropy provides valuable insights into the complexity and diversity of life, simplifying the study for researchers to study and understand biological systems. References \citep{jost2006entropybio1}, \citep{carranza2007analyzingbio2}, \citep{gao2019computationbio3}, and \citep{gnaiger1989physiologicalbio4} provide a detailed discussion of the applicability of entropy measures in bioenergetics, ecology, and evolutionary biology, as discussed above.

Entropy measures also find applications in fields such as linguistics, psychology, sociology, communication studies, and visual and performing arts. However, since their usage is mostly limited to measuring randomness and uncertainty, we are not including related literature on these applications.

\section{Data Resources for Entropy Measure Applications}\label{section4}
This section conducts a thorough analysis of the recent and prevalent applications of entropy measures across diverse datasets. We have meticulously examined each reference in this article regarding the application of entropy to datasets. This evaluation led us to identify $81$ references that apply entropy measures. Furthermore, we delve deeper into our analysis to focus specifically on openly accessible datasets, yielding a total of $49$ references corresponding to datasets available through open sources or directories. On the other hand, we categorize private datasets according to their availability; these include datasets that either the respective authors have not shared or are not accessible through the indicated sources. Figure \ref{fig:8} illustrates the approximate proportion of data attributed to openly accessible datasets in comparison to that of private datasets.

\begin{figure}
    \centering
    \includegraphics[width=0.7\linewidth]{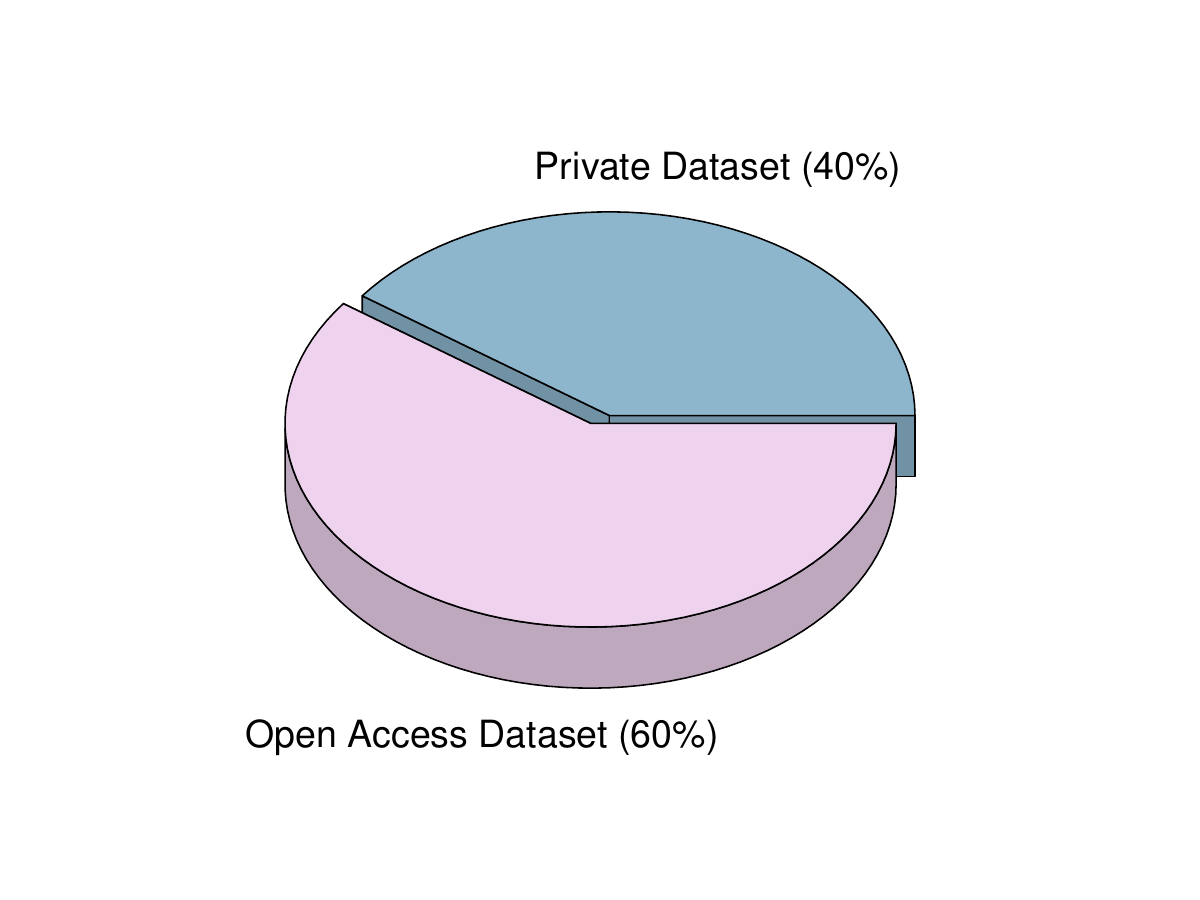}
    \caption{Sharing of Open Access and Private datasets for all entropy application references}
    \label{fig:8}
\end{figure}

Table \ref{Tableofdatarepositories} presents an overview of references employing various entropy measures in their applications, detailing the types of entropy, sources, and links to access the openly available datasets. The following list of $32$ references includes private datasets that could be useful for their methodology:

\citep{al2020novelIMG34}, \citep{al2018newIMG55}, \citep{baghmolaei2022applicationCRE1finance3}, \citep{tranferentropybook}, \citep{cao2017inherentfuzzy2Time1}, \citep{carranza2007analyzingbio2}, \citep{chen2001musicMRS}, \citep{chen2007characterizationfuzzy1Time1}, \citep{das2017detectingshannoncategorical3}, \citep{de2004imageIMG22},
\citep{dimpfl2013usingTransferentropyfinance3}, \citep{fadlallah2013weightedpermutationentropy}, \citep{fraiwan2021gaugingT2AI1}, \citep{hasan2010regisrationIMG32}, \citep{jost2006entropybio1}, \citep{kaniadakis2002statisticalkaniadakisentropy}, \citep{nasr2002cross3AI1}, \citep{rahmani2020kullbackleiblerdivergencecommunication2}, \citep{rao2004cumulativeresidualentropy}, \citep{raychaudhury1984discriminationRRGRBgraphentropy}, \citep{richman2000physiologicalsampleentropydefination}, \citep{romero2001sequencecomplexityAI1}, \citep{sabeti2009entropycomplexitymeasureAI1}, \citep{sahoo2004thresholdingIMG21}, \citep{shi2013semgresidual1Time1}, \citep{shiozaki1986edgeIMG12}, \citep{simovici2000generalizedshannoncategorical5}, \citep{sun2022statisticalCRE4finance3}, \citep{valencia2009refinedmultiscaleentropy}, \citep{wang1984analysisdecisiontreeAI2}, \citep{zhan2015medicalIMG52}, \citep{zhan2016graphIMG51}.

\begin{figure}
    \centering
    \includegraphics[width=1\linewidth]{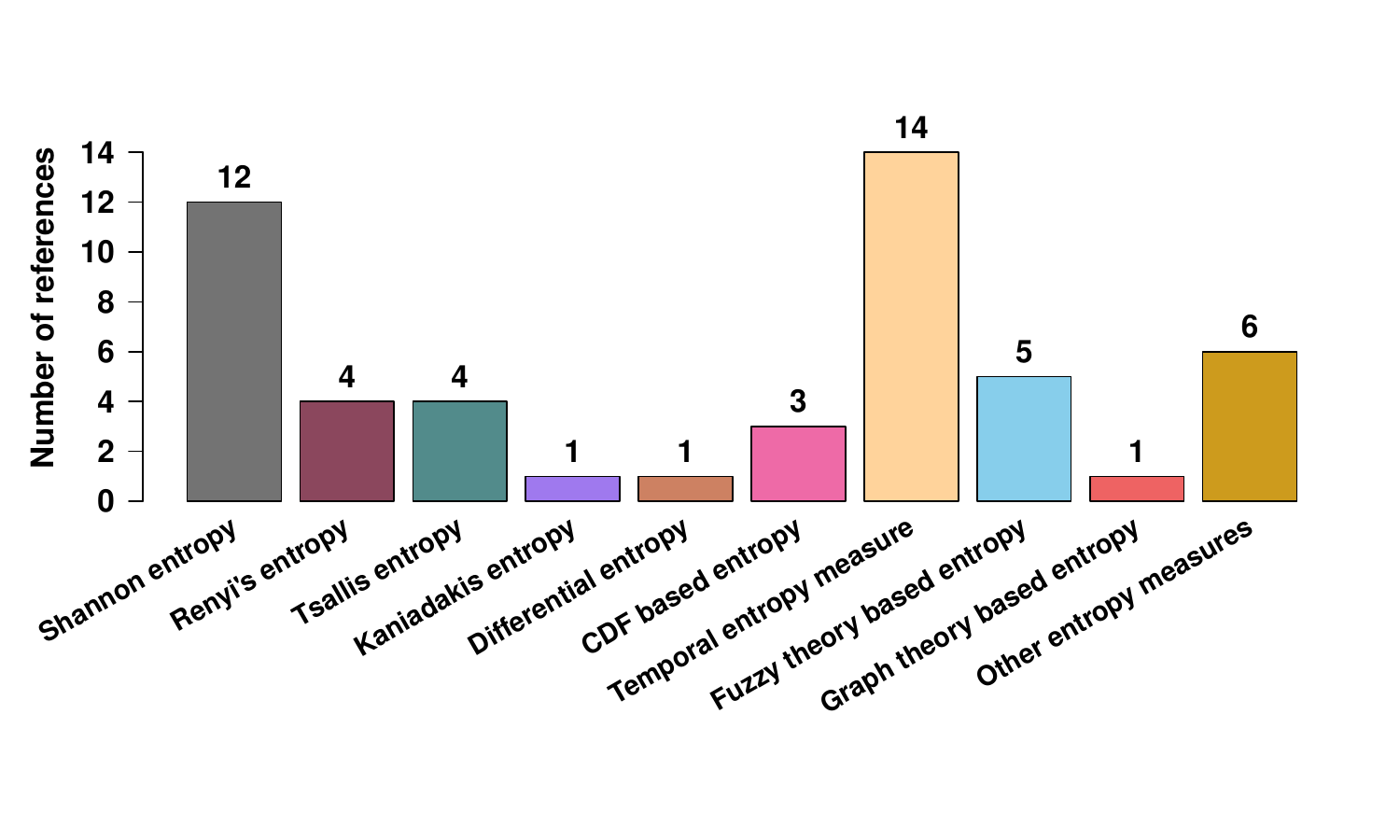}
    \caption{Spread of various entropy measures in openly available datasets}
    \label{fig:9}
\end{figure}

\begin{figure}[H]
    \centering
    \includegraphics[width=1\linewidth]{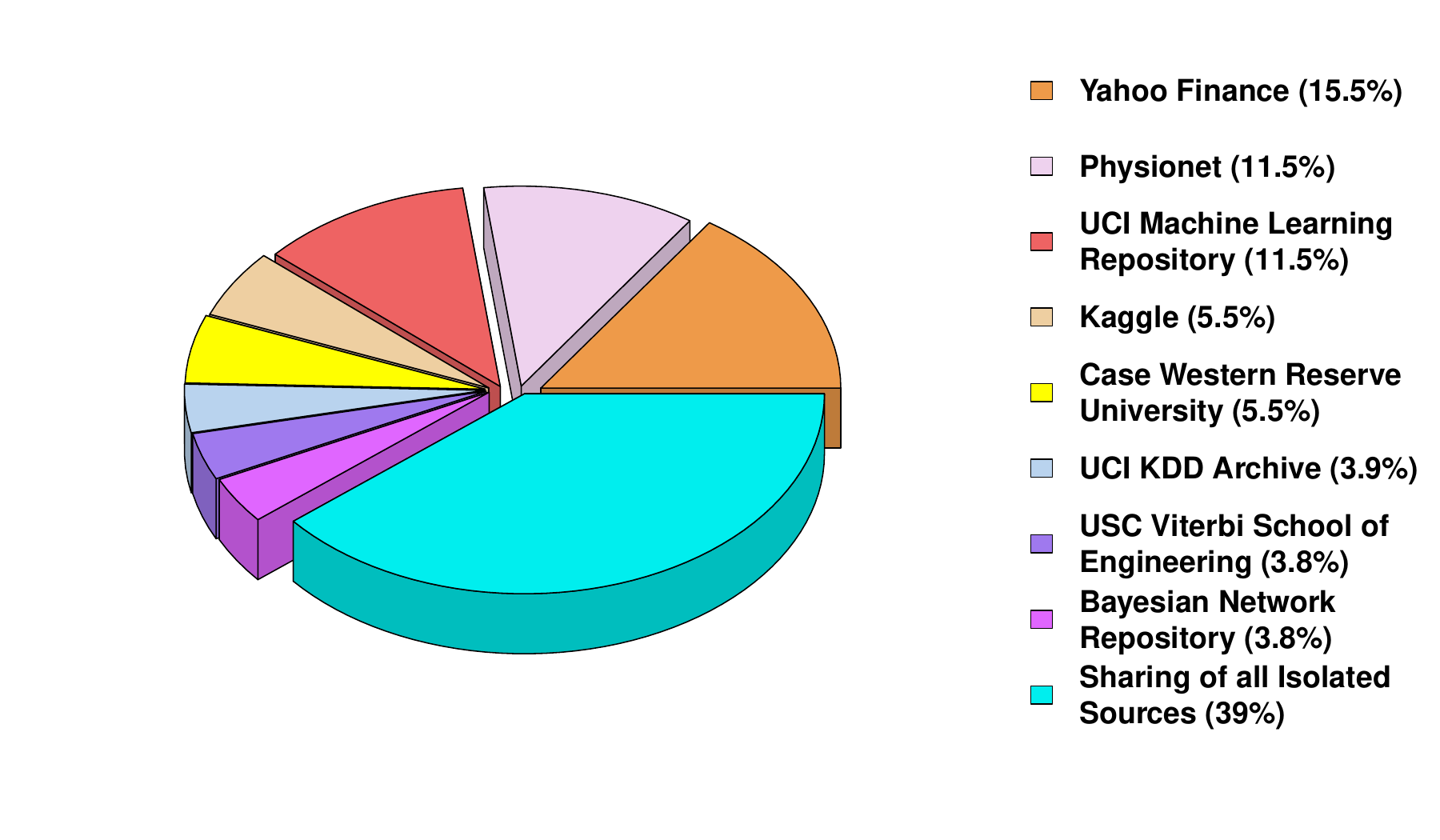}
    \caption{Allocation of various data source repositories across all openly available dataset references}
    \label{fig:10}
\end{figure}

Furthermore, we analyze various types of entropy measures applied to openly accessible datasets to see recent trends in present applications. Figure \ref{fig:9} presents the corresponding analysis, highlighting the importance of Temporal Entropy Measures and Shannon Entropy. These insights are derived from Table \ref{Tableofdatarepositories}, which categorizes the specific entropy measures used in the corresponding references. We categorize the sources of openly available datasets, with Figure \ref{fig:10} showing the proportional contribution of each. The 'Sharing of all Isolated Reference Sources' includes those from Table \ref{Tableofdatarepositories} with only a single reference. Notably, 'Yahoo Finance,' 'Physionet,' and the 'UCI Machine Learning Repository' together account for $38.5\%$ of the total dataset.

\begin{sidewaystable}
    \centering
    \begin{tabular}{p{6cm} p{9cm} p{7.5cm}} 
        \toprule
        \textbf{Reference} & \textbf{In-use entropy} & \textbf{Source link}\\
        \midrule
        \citep{ahmed2011multivariatemultiscalentropy},\citep{alcaraz2010reviewsampleTime1},\citep{aziz2005multiscalepermutationentropy},\citep{costa2002multiscaleentropydefination1},\citep{costa2005multiscaleentropydefination2},\citep{fabila2022permutationgraph2Time1} & Temporal entropy measure  & \href{https://physionet.org/about/database/}{\textcolor{blue}{Physionet}}\\
        \cmidrule(lr){1-3}
        \citep{alsakran2014usingshannonentropycategorical2},\citep{cui2019improveddengAI1},\citep{mahamadou2020categorical7},\citep{morente2016improvingfuzzy1AI1},\citep{song2022combininggraph1Time1},\citep{szmidt2014measureintuitionisticfuzzyentropy} & Shannon, Deng, Fuzzy and Belief entropy & \href{https://archive.ics.uci.edu/datasets}{\textcolor{blue}{UCI Machine Learning Repository}}\\
        \cmidrule(lr){1-3}
        \citep{machadorenyifractionalorderentropy},\citep{machadofractionalorderentropy},\citep{olbrys2022approximatesampleentropyfinance3},\citep{wang2021generalizeddispersionentropyfinance3},\citep{xia2012multiscaleentropyfinance3},\citep{xiong2019fractionalcumulativeresidualentropy},\citep{yang2012studyCRE2finance3},\citep{yin2016weightedpermutationentropyproperties} & Fractional order, Temporal, Cumulative residual entropy & \href{https://finance.yahoo.com/}{\textcolor{blue}{Yahoo finance}}\\
        \cmidrule(lr){1-3}

        \citep{barbara2002coolcatcategorical6},\citep{lima2010decisiontreeAI1} & Shannon, Renyi's and Tsallis entropy & \href{https://kdd.ics.uci.edu/summary.data.alphabetical.html}{\textcolor{blue}{UCI KDD archive}}\\
        \cmidrule(lr){1-3}
        \citep{ho2019realcross2AI1},\citep{lima2010decisiontreeAI1},\citep{kalra2010automaticTsallis1Time1} & Shannon, Renyi's, Tsallis and Cross entropy  & \href{https://www.kaggle.com/}{\textcolor{blue}{Kaggle}}\\
        \cmidrule(lr){1-3}
        \citep{sharma2016feature2AI1},\citep{wu2013modifiedmultiscaleentropy},\citep{wu2013timecompositemultiscaleentropy} & Shannon and Temporal entropy & \href{https://engineering.case.edu/bearingdatacenter/download-data-file}{\textcolor{blue}{Case Western Reserve University Ohio}}\\
        \cmidrule(lr){1-3}

        \citep{zhang2015pathologicalTsallis2AI1} & Tsalli entropy & \href{https://dataverse.harvard.edu/}{\textcolor{blue}{Harvard Dataverse}}\\
        \cmidrule(lr){1-3}

        \citep{beenamol2012waveletTsallis2Time1} & Shannon and Tsallis entropy & \href{https://www.strongmotioncenter.org/index.html}{\textcolor{blue}{Center for Engineering Strong Motion Data}}\\
        \cmidrule(lr){1-3}

        \citep{xiang2020imageIMG33} & Residual entropy & \href{https://yann.lecun.com/exdb/mnist/}{\textcolor{blue}{Yann LeCun's website}}\\
        \cmidrule(lr){1-3}

        \citep{sanyal2011adaptiveIMG42},\citep{wu2013localIMG13} & Shannon and Fuzzy entropy & \href{https://sipi.usc.edu/database/}{\textcolor{blue}{USC Viterbi School of Engineering}}\\
        \cmidrule(lr){1-3}

        \citep{rostaghi2019applicationdispersionTime1} & Temporal entropy & \href{https://cir.nii.ac.jp/crid/1571980075951352192}{\textcolor{blue}{CiNii research}}\\
        \cmidrule(lr){1-3}

        \citep{compton2022entropicCI3},\citep{compton2020entropicCI2} & Shannon entropy & \href{https://www.bnlearn.com/bnrepository/}{\textcolor{blue}{Bayesian Network Repository}}\\
        \cmidrule(lr){1-3}

        \citep{dehmer2011historygraphentropychemistry1} & Graph entropy & \href{https://data.nist.gov/sdp/#/}{\textcolor{blue}{NIST Data}}\\
        \cmidrule(lr){1-3}

        \citep{fabila2022permutationgraph2Time1} & Temporal entropy & \href{https://www.data.gouv.fr/fr/datasets/}{\textcolor{blue}{Open data}}\\
        \cmidrule(lr){1-3}

        \citep{guido2018tutorialfeature1AI1} & Shannon entropy & \href{https://www.ibilce.unesp.br/#!/departamentos/cienc-comp-estatistica/docentes/rodrigo-capobianco-guido/}{\textcolor{blue}{Dr. Eng. Rodrigo Capobianco Guido}}\\
        \cmidrule(lr){1-3}

        \citep{huang2022permutationentropyproperties2} & Temporal entropy & \href{https://www.lseg.com/en/data-analytics/financial-data}{\textcolor{blue}{LSEG data and analytics}}\\    
        \cmidrule(lr){1-3}
        
        \citep{lei2020adaptiveIMG24} & Kaniadakis entropy & \href{http://vcipl-okstate.org/pbvs/bench/}{\textcolor{blue}{VCIPL of Oklahoma State University}}\\

         \cmidrule(lr){1-3}
        
          \citep{machadorenyifractionalorderentropy} & Fractional entropy & \href{https://www.scopus.com/search/form.uri?display=basic#basic}{\textcolor{blue}{Scopus database}} \\

         \cmidrule(lr){1-3}
          \citep{kalra2010automaticTsallis1Time1} & Tsallis entropy & \href{https://www.library.ucsf.edu/}{\textcolor{blue}{UCSF library}} \\

         \cmidrule(lr){1-3}

          \citep{naidu2018shannonIMG44} & Shannon and Fuzzy entropy & \href{https://www.imageprocessingplace.com/root_files_V3/image_databases.htm}{\textcolor{blue}{Pearson-Prentice Hall image database}} \\

         \cmidrule(lr){1-3}

          \citep{rostaghi2016dispersionentropy} & Temporal entropy & \href{https://www.upf.edu/web/mdm-dtic/datasets-2016}{\textcolor{blue}{Universitat Pompeu Fabra, Barcelona}} \\

         \cmidrule(lr){1-3}

          \citep{sanyal2011adaptiveIMG42}  & Fuzzy entropy & \href{https://ccia.ugr.es/cvg/CG/base.htm}{\textcolor{blue}{DCSAI, University of Granada}} \\

         \cmidrule(lr){1-3}

          \citep{arroyo2019unsupervisedsentencerepresentationAI1} & Shannon entropy & \href{https://wiki.cimec.unitn.it/tiki-index.php?page=CLIC}{\textcolor{blue}{CLIC wikipage}} \\

         \cmidrule(lr){1-3}

          \citep{arroyo2019unsupervisedsentencerepresentationAI1} & Shannon entropy & \href{https://ixa2.si.ehu.eus/stswiki/}{\textcolor{blue}{STS benchmark dataset}} \\

         \cmidrule(lr){1-3}
          \citep{silva2016twoIMG53} & Temporal entropy & \href{https://multibandtexture.recherche.usherbrooke.ca/original_brodatz.html}{\textcolor{blue}{Université de Sherbrooke Canada}}\\

         \cmidrule(lr){1-3}
          \citep{wang2007nonIMG31} & Cumulative residual entropy & \href{https://brainweb.bic.mni.mcgill.ca/}{\textcolor{blue}{Brainweb}} \\

         \cmidrule(lr){1-3}
          \citep{xu2021deepgraph1AI1} & Shannon and Renyi's entropy & \href{https://huggingface.co/datasets/graphs-datasets/MUTAG}{\textcolor{blue}{Graph dataset from hugging face}}\\

         \cmidrule(lr){1-3}
          \citep{xu2021deepgraph1AI1} & Shannon and Renyi's entropy & \href{https://huggingface.co/datasets/graphs-datasets/MUTAG}{\textcolor{blue}{Network Repository}}\\

        \bottomrule
    \end{tabular}
    \caption{A list of references of all studies based on openly accessible datasets, detailing the applied entropy measures and providing links to the corresponding datasets.}
    \label{Tableofdatarepositories}
\end{sidewaystable}

\section{Future Trends}\label{section5}
The primary goal of entropy is to capture the uncertainty, randomness, or complexity inherent in a system, thereby aiding in the modelling of its behaviour. As a result, many mathematical measures have been developed based on various influencing factors. Below are some potential areas and associated challenges where entropy might play a significant role:
\begin{enumerate}
    \item Recently, high-entropy alloys\citep{george2019definationofhighentropyalloysthermodynamics1} have gained significant attention due to their numerous benefits. In addition to enhancing strength and flexibility, other aspects such as corrosion resistance, oxidation resistance, magnetic properties, and material degradation also require further exploration. The configurational entropy can be tailored to design materials with desired characteristics for specific applications.

    \item OpenAI has experienced rapid growth in a short period, with its applications expanding across numerous technical fields to enhance efficiency. However, challenges still need to be addressed, including outdated information, inconsistent outputs, difficulty in handling multimodal tasks, the black-box nature of its models, and concerns over plagiarism and copyright issues. The role of entropy in addressing these challenges has yet to be fully explored, but it may offer a promising approach in managing such issues more effectively. 

    \item In every definition of entropy, the distribution function corresponding to a natural phenomenon is used to compute entropy values. This approach often overlooks the significance of the nature of the random variable itself, which can offer a new perspective on disorder or complexity. Therefore, generalizing the definition of entropy to account for the characteristics of the underlying random variable could be a valuable direction for modelling more realistic phenomena.

    \item The digital revolution is rapidly progressing worldwide, as evidenced by movies, OTT content, TV shows, and tools such as search engines, cloud storage services, productivity tools, and social media platforms. A common challenge across these services is to provide accurate and personalized recommendations\citep{chen2001musicMRS}. Entropy can be leveraged to capture individual user complexity, enhancing the effectiveness of recommended systems, and delivering more tailored content.

    \item The Principle of Maximum Entropy is a widely used method for identifying the most probable state of a system or the distribution of data by utilizing various entropy measures under given constraints. However, a key limitation of this approach is that it may not provide a precise approximation of reality without accurate knowledge of the constraints. To improve the realism of system distribution estimation, new techniques, such as AI-driven models, for feeding accurate data information and others, can be explored and applied across various domains.

    \item In causal inference problems, entropy can help in determining the direction of causality by considering factors such as conditional independence of random variables, feature selection, and treatment effects. Despite significant efforts over the past five years\citep{kocaoglu2020applicationsCI1,compton2020entropicCI2,compton2022entropicCI3}, exploring entropy's potential in causal inference is still in its early stages.
    
    \item Duality in mathematics provides alternative perspectives on problems, uncovering hidden symmetries, simplifying solutions, and offering deeper theoretical insights. Extropy\citep{lad2015extropy}, introduced in 2015 as the complementary dual of entropy, measures underlying uncertainty in systems. However, the differences between various forms of entropy and extropy, particularly in terms of their physical interpretations and relationships, remain largely unexplored.
    
\end{enumerate}

\section{Conclusion}\label{section6}
The present review addresses the need for a concise view of entropy, providing a foundation for researchers to explore directions that can contribute to the scientific community for real-world applications, and the identification of new research interests. The article provides an in-depth study of $60$ entropy measures, with respective motivations, fundamental properties, and interrelations. Although the total number of entropy measures in the literature is far larger, we selected the specific ones based on their prominence and wide applications, as indicated by the citation counts. These measures are categorized into seven categories based on their underlying computational theory. It is noteworthy that concepts from interdisciplinary paradigms have driven the significant development of entropy measures over the past $20$ years. Additionally, the article presents contemporary applications across eight major fields that are currently benefited from entropy. The application section effectively conveys that entropy can be utilized in real-life erratic phenomena. A collection of $81$ references is provided, listing the use of entropy measures directly applied to datasets, with $49$ openly accessible datasets and their resource links. The article also uncovers meaningful new directions for researcher in entropy and thus, may provide a guided path to the researchers in the field.

\bibliography{Refreview}

\end{document}